\documentclass[11pt]{article}
\usepackage[english]{babel}
\usepackage{amsmath,amssymb,latexsym,theorem,graphicx,enumerate}
\usepackage[margin=2cm]{geometry}
\newcommand{\assign}{:=}
\newcommand{\backassign}{=:}
\newcommand{\cdummy}{\cdot}
\newcommand{\longhookrightarrow}{{\lhook\joinrel\relbar\joinrel\rightarrow}}
\newcommand{\nin}{\not\in}
\newcommand{\of}{:}
\newcommand{\suchthat}{:}

\newcommand{\tmem}[1]{\emph{#1}}

\newcommand{\tmname}[1]{\textsc{#1}}
\newcommand{\tmnote}[1]{\thanks{\ #1}}
\newcommand{\tmop}[1]{\ensuremath{\operatorname{#1}}}
\newcommand{\tmrsub}[1]{\ensuremath{_{\textrm{#1}}}}
\newcommand{\tmsamp}[1]{\textsf{#1}}
\newcommand{\tmtextit}[1]{\textit{#1}}
\newcommand{\tmstrong}[1]{\textbf{#1}}
\newcommand{\tmtextup}[1]{\text{{\upshape{#1}}}}
\newenvironment{descriptioncompact}{\begin{description} }{\end{description}}
\newenvironment{enumerateroman}{\begin{enumerate}[i.] }{\end{enumerate}}
\newenvironment{itemizedot}{\begin{itemize} }{\end{itemize}}
\newenvironment{proof}{\noindent\textbf{Proof\ }}{\hspace*{\fill}$\Box$\medskip}
\newenvironment{tmcompact}{\begin{tmparsep}{0em}}{\end{tmparsep}}
\newenvironment{tmparmod}[3]{\begin{list}{}{\setlength{\topsep}{0pt}\setlength{\leftmargin}{#1}\setlength{\rightmargin}{#2}\setlength{\parindent}{#3}\setlength{\listparindent}{\parindent}\setlength{\itemindent}{\parindent}\setlength{\parsep}{\parskip}} \item[]}{\end{list}}
\newenvironment{tmparsep}[1]{\begingroup\setlength{\parskip}{#1}}{\endgroup}
\newcounter{tmcounter}
\newcommand{\custombinding}[1]{%
  \setcounter{tmcounter}{#1}%
  \addtocounter{tmcounter}{-1}%
  \refstepcounter{tmcounter}}
{\theorembodyfont{\rmfamily}\newtheorem{example}{Example}}
\newtheorem{lemma}{Lemma}
\newtheorem{proposition}{Proposition}
{\theorembodyfont{\rmfamily}\newtheorem{remark}{Remark}}
\newtheorem{theorem}{Theorem}
\newcommand{\tmkeywords}{\textbf{Keywords:} }
\usepackage{url}

\begin{document}

\title{A proof of the Brill--Noether method from scratch%
  \tmnote{This paper is part of a project that has received funding from the
  French ``Agence de l'innovation de d{\'e}fense''. Elena Berardini has
  received funding from the European Union's Horizon 2020 research and
  innovation programme under the Marie Sk{\l}odowska-Curie grant agreement No
  899987. Alain Couvreur is funded by the ANR Grant
  ANR-21-CE39-0009-BARRACUDA.}
}

\author{Elena Berardini\smallskip\\
CNRS; Institut de Mathématiques de Bordeaux\\
Universit{\'e} de Bordeaux\\
  Talence, France, F 33405 \smallskip\\
Department of Mathematics and Computer Science\\
  Eindhoven University of Technology\\
  Eindhoven, the Netherlands, 5612 AZ\\
  \url{elena.berardini@math.u-bordeaux.fr}
  \and
Alain Couvreur\smallskip \\
Inria \& Laboratoire d'informatique de l'{\'E}cole polytechnique (LIX, UMR 7161)\\
  CNRS, {\'E}cole polytechnique, Institut Polytechnique de Paris\\
  Palaiseau, France, 91120\\
  \url{alain.couvreur@inria.fr} 
  \and
  Gr{\'e}goire Lecerf\smallskip\\
  Laboratoire d'informatique de l'{\'E}cole polytechnique (LIX, UMR 7161)\\
  CNRS, {\'E}cole polytechnique, Institut Polytechnique de Paris\\
  Palaiseau, France, 91120\\
  \url{gregoire.lecerf@lix.polytechnique.fr}
  }

\date{}

\maketitle

\begin{abstract}
  In 1874 Brill and Noether designed a seminal geometric method for computing
  bases of Riemann--Roch spaces. From then, their method has led to several
  algorithms, some of them being implemented in computer algebra systems. The
  usual proofs often rely on abstract concepts of algebraic geometry and
  commutative algebra. In this paper we present a short self-contained and
  elementary proof that mostly needs Newton polygons, Hensel lifting,
  bivariate resultants, and Chinese remaindering.
\end{abstract}

\tmkeywords{Algebraic curves, Riemann--Roch spaces, Brill--Noether method,
Hensel lemmas, Newton polygons}

\section{Introduction}

Riemann--Roch spaces are vector spaces of rational functions that satisfy some
conditions on the localization and the multiplicity of their zeros and poles.
These spaces are a cornerstone of modern applications of algebra to various
areas of computer science. This paper presents the seminal method due to Brill
and Noether to compute Riemann--Roch spaces. Our approach here is mostly
dedicated to undergraduate students: it relies on Newton polygons, Hensel
lifting, bivariate resultants, and Chinese remaindering.

This introduction begins with elementary problems in order to motivate the
study of Riemann--Roch spaces and to sketch the Brill--Noether method.
Section~\ref{s:preliminaries} gathers well known elementary definitions and
results from algebra, while Section~\ref{s:hensel} is dedicated to the Hensel
lemma and its usual extensions. These tools are central in studying how
valuations extend to field extensions and to introduce the notion of
uniformizing parameter, that is the goal of Section~\ref{s:valuations}. All
these results allow us to define places and divisors on curves in
Section~\ref{s:divisor}. Finally Section~\ref{s:BNmethod} is devoted to the
Brill--Noether algorithm. The proofs are self-contained. References to
algorithms in the literature are gathered in our final
Section~\ref{s:notes}.{\smallskip}

Let $\mathbb{K}$ be a field. The ring of univariate polynomials over
$\mathbb{K}$ is written $\mathbb{K} [x]$, and we let $\mathbb{K}
[x]_{\leqslant d} \assign \{ f \in \mathbb{K} [x] : \deg f \leqslant d \}$.
The ring of {\tmem{power series}} over $\mathbb{K}$ is written $\mathbb{K}
[[x]]$ and its field of fractions, called the field of {\tmem{Laurent
series}}, is denoted by $\mathbb{K} ((x))$.

Let $\mathbb{D}$ be an integral domain. A {\tmem{valuation}} $v$ on
$\mathbb{D}$ is a map from $\mathbb{D}$ to $\Gamma \cup \{ \infty \}$ such
that:

\begin{tmcompact}
  \begin{itemizedot}
    \item $\Gamma$ is an abelian additive group endowed with a total order
    written $<$;
    
    \item $v (a) = \infty$ if and only if $a = 0$;
    
    \item $v (ab) = v (a) + v (b)$;
    
    \item $v (a + b) \geqslant \min (v (a), v (b))$.
  \end{itemizedot}
\end{tmcompact}

The field $\mathbb{K} ((x))$ is endowed with the natural valuation in $x$,
written $\tmop{val}_x$, normalized such that $\tmop{val}_x (x) = 1$, so
$\tmop{val}_x \left( \sum_{i \in \mathbb{Z}} a_i x^i \right) \assign \min (i
\in \mathbb{Z} \suchthat a_i \neq 0)$.

\subsection{Riemann--Roch spaces}\label{s:intro-rr}

Let $\alpha_1, \ldots, \alpha_r$ be distinct values in $\mathbb{K}$. As a
basic fact of algebra, it is well known that the set of polynomials $g \in
\mathbb{K} [x]$ that vanish simultaneously at $\alpha_1, \ldots, \alpha_r$
with respective multiplicities $m_1, \ldots, m_r$ is the ideal of $\mathbb{K}
[x]$ generated by $(x - \alpha_1)^{m_1} \cdots (x - \alpha_r)^{m_r}$. The
polynomials in this ideal that satisfy a given degree upper bound constitute a
finite dimensional $\mathbb{K}$\mbox{-}vector space, of which a basis is
easily given. This setting extends to rational functions, that is, fractions
of polynomials, and signed integers $m_i$. In this context, we can pose the
following problem.

\begin{descriptioncompact}
  \item[Problem A] Let $\alpha_1, \ldots, \alpha_r$ be distinct values in
  $\mathbb{K}$, let $m_1, \ldots, m_r$ be signed integers, and set $D_A
  \assign \{ (\alpha_1, m_1), \ldots, (\alpha_r, m_r) \}$. Compute the
  $\mathbb{K}$-vector space $\mathcal{L} (D_A)$ of rational functions $g / h
  \in \mathbb{K} (x)$ such that $\deg g \leqslant \deg h$, the valuation of $g
  / h$ regarded in $\mathbb{K} ((x - \alpha_i))$ is~$\geqslant - m_i$ for $i =
  1, \ldots, r$, and that for any other value $\alpha$ in $\mathbb{K}$ the
  valuation of $g / h$ regarded in $\mathbb{K} ((x - \alpha ))$ is~$\geqslant
  0$.
\end{descriptioncompact}

It is clear that $\mathcal{L} (D_A)$ is a $\mathbb{K}$-vector space. A basis
of $\mathcal{L} (D_A)$ can be computed as follows. We let $d \assign \sum_{i =
1, m_i > 0}^r m_i$ and $h (x) \assign \prod_{i = 1, m_i > 0}^r (x -
\alpha_i)^{m_i}$. If $a / b \in \mathcal{L} (D_A)$ with $a$ and~$b$ coprime,
then $b$ divides $h$, hence there exists a unique polynomial $g$ of
degree~$\leqslant d$ such that $g / h = a / b$. This shows that $\mathcal{L}
(D_A)$ has finite dimension and that $h$ can be taken as a common denominator
of a basis of $\mathcal{L} (D_A)$. So we are led to seek for a basis $g_1,
\ldots, g_{\ell}$ of polynomials $g \in \mathbb{K} [x]_{\leqslant d}$ such
that the valuation of $g$ regarded in $\mathbb{K} [[x - \alpha_i]]$ is
$\geqslant - m_i$ for $i = 1, \ldots, r$. This task reduces to solving a
linear system in the $d + 1$ unknown coefficients of~$g$. If the number of
equations, namely $\sum_{i = 1, m_i < 0}^r (- m_i)$, does not exceed the
number of unknowns, then the set of non-zero solutions is not empty. More
precisely we may take
\[ g_1 (x) \assign \prod_{i = 1, m_i < 0}^r (x - \alpha_i)^{- m_i} \text{ and
   } g_i \assign x^{i - 1} g_1 \text{ for } i = 2, \ldots, \ell, \text{ where}
\]
\begin{equation}
  \ell \assign (d + 1) - \sum_{i = 1, m_i < 0}^r (- m_i) = 1 + \sum_{i = 1}^r
  m_i . \label{eq:intro-ell}
\end{equation}
Finally, $g_1 / h, \ldots, g_{\ell} / h$ is a basis of $\mathcal{L} (D_A)$.

Problem~A naturally extends to the projective setting: the affine line is
replaced by the projective line~$\mathbb{P}^1 (\mathbb{K})$ over $\mathbb{K}$.
Rational functions defined on $\mathbb{P}^1 (\mathbb{K})$ are either $0$ or of
the form $A / B$ where $A$ and $B$ are homogeneous polynomials in $\mathbb{K}
[x, z]$ of same degree. The set of such functions is written $\mathbb{K}
(\mathbb{P}^1 (\mathbb{K}))$. If $\zeta_i = (\alpha_i \of 1)$ is a point in
the affine chart $z = 1$ of $\mathbb{P}^1 (\mathbb{K})$, then the valuation of
$G / H$ at~$\zeta_i$ is the valuation of $G (x, 1) / H (x, 1)$ regarded in
$\mathbb{K} ((x - \alpha_i))$. If $\zeta_i$ is the point $(1 : 0)$, called the
point ``at infinity'', then the valuation $G / H$ at $\zeta_i$ is the
valuation of $G (1, z) / H (1, z)$ regarded in~$\mathbb{K} ((z))$.

\begin{descriptioncompact}
  \item[Problem B] Let $\zeta_1, \ldots, \zeta_r$ be distinct points on the
  projective line $\mathbb{P}^1 (\mathbb{K})$, let $m_1, \ldots, m_r$ be
  signed integers, and set $D_B \assign \{ (\zeta_1, m_1), \ldots, (\zeta_r,
  m_r) \}$. Compute the $\mathbb{K}$-vector space $\mathcal{L} (D_B)$ of
  rational functions $G / H \in \mathbb{K} (\mathbb{P}^1 (\mathbb{K}))$ such
  that the valuation of $G / H$ at $\zeta_i$ is~$\geqslant - m_i$ for $i = 1,
  \ldots, r$, and that for any other point $\zeta$ in $\mathbb{P}^1
  (\mathbb{K})$ the valuation of $G / H$ at $\zeta$ is $\geqslant 0$.
\end{descriptioncompact}

If all the $\zeta_i$ are in the affine chart $z = 1$, then we may write them
as $\zeta_i = (\alpha_i \suchthat 1)$ and reduce Problem~B to Problem~A: let
$g_1 / h, \ldots, g_{\ell} / h$ be a basis of $\mathcal{L} (D_A)$, then with
$d \assign \deg h$ and $H \assign z^d h (x / z)$, we verify routinely that
$z^d g_1 (x / z) / H, \ldots, z^d g_{\ell} (x / z) / H$ is a basis of
$\mathcal{L} (D_B)$. If one of the $\zeta_i$ is $(1 : 0)$ and if the
cardinality of $\mathbb{K}$ is sufficiently large, then we may apply a
sufficiently generic linear change of variables
\[ P_i = (\gamma_i \of \delta_i) \mapsto (a \gamma_i + b \delta_i \of c
   \gamma_i + d \delta_i) = \left( \frac{a \gamma_i + b \delta_i}{c \gamma_i +
   d \delta_i} \of 1 \right) \]
with $ad - bc \neq 0$ and such that $c \gamma_i + d \delta_i \neq 0$ for all
$i = 1, \ldots, r$, in order to reduce again Problem~B to Problem~A.

\begin{remark}
  If one reformulates Problem~A into the projective setting, the condition
  $\deg g \leqslant \deg h$ is equivalent to having a nonnegative valuation at
  infinity. Thus, Problem~B can be regarded as a generalization of Problem~A
  where any condition at infinity can be imposed.
\end{remark}

Alternatively, we may solve Problem~B directly as follows. We still set $d
\assign \sum_{i = 1, m_i > 0}^r m_i$ but set $H (x, z) \assign \prod_{i = 1,
m_i > 0}^r (\beta_i x - \alpha_i z)^{m_i}$, where $\alpha_i$ and $\beta_i$
represent coordinates of~$\zeta_i = (\alpha_i \of \beta_i)$. If $A / B \in
\mathcal{L} (D_B)$ with $A$ and $B$ coprime, then $B$ divides $H$ and there
exists a unique homogeneous polynomial $G$ of degree $d$ such that $G / H = A
/ B$. This shows that $\mathcal{L} (D_B)$ has finite dimension and that $H$
can be taken as a common denominator of a basis of~$\mathcal{L} (D_B)$. So it
remains to seek for a basis $G_1, \ldots, G_{\ell}$ in $\mathbb{K} [x, y, z]$
of homogeneous polynomials~$G$ of degree~$d$ such that the valuation of $G /
H$ at $P_i$ is $\geqslant - m_i$ for $i = {1, \ldots, r}$. This task reduces
to solving a linear system in the $d + 1$ unknown coefficients of $G$. Finally
$G_1 / H, \ldots, G_{\ell} / H$ is a basis of $\mathcal{L} (D_B)$, where
$\ell$ is the same as in~\eqref{eq:intro-ell}.

\begin{example}
  \label{ex:rr-rs}A particularly simple case of interest of Problem~B is the
  following: $r = 1$, $\zeta_1 = (1 : 0)$ and $m_1 \geqslant 0$. Then, the
  above direct method shows that $x^{m_1} / z^{m_1}, x^{m_1 - 1} / z^{m_1 -
  1}, \ldots, x / z, 1$ is a basis of~$\mathcal{L} (D_B)$.
\end{example}

Problem~B can be further generalized to plane algebraic curves. Let $F \in
\mathbb{K} [x, y, z]$ be an irreducible homogeneous polynomial. Let
$\mathbb{P}^2 (\mathbb{K})$ stands for the projective plane over~$\mathbb{K}$.
The set of zeros of $F$ in $\mathbb{P}^2 (\mathbb{K})$ is written
$\mathcal{C}$; this is an {\tmem{algebraic projective plane curve}}. A
point~$\zeta$ of $\mathcal{C} \mathcal{}$ is called \tmtextit{singular} if the
partial derivatives of $F$ vanish simultaneously at $\zeta$, and
\tmtextit{non-singular} otherwise. The set of singular points of a curve is
called its \tmtextit{singular locus.} The set of rational functions defined on
$\mathcal{C}$ over $\mathbb{K}$, written $\mathbb{K} (\mathcal{C})$, is the
set of functions that are either $0$ or of the form $A / B$ where $A$ and $B$
are homogeneous polynomials of the same degree, with $B$ prime to~$F$, and
subject to the equivalence relation
\[ A / B \sim A' / B' \Longleftrightarrow AB' - A' B \in (F) . \]
At a non-singular point $\zeta = (\zeta_x : \zeta_y : \zeta_z)$ of
$\mathcal{C}$, and up to an occasional linear change of variables, we may
assume that $\zeta_z = 1$ and $\frac{\partial F}{\partial y} (\zeta) \neq 0$,
so the implicit function theorem ensures that $\mathcal{C}$ is locally defined
by a power series
\[ \varphi (x) \assign \zeta_y + c_1  (x - \zeta_x) + c_2  (x - \zeta_x)^2 +
   \cdots \in \mathbb{K} [[x - \zeta_x]] \]
that satisfies $F (x, \varphi (x), 1) = 0$; see~{\cite[Chapter~7, Section~7.2,
Corollary~7.4]{Eisenbud1995}}. If $G / H \in \mathbb{K} (\mathcal{C})$, then
its valuation at $\zeta$ is defined by
\[ \tmop{val}_{\zeta} (G / H) \assign \tmop{val}_{x - \zeta_x} (G (x, \varphi
   (x), 1) / H (x, \varphi (x), 1)) ; \]
we will see that this definition is essentially independent of the choice of
coordinates. Defining valuations at singular points is more intricate: this
will be the purpose of Section~\ref{s:divisor}.

\begin{descriptioncompact}
  \item[Problem C] Let $\zeta_1, \ldots, \zeta_r$ be distinct non-singular
  points on the projective plane curve $\mathcal{C}$, let $m_1, \ldots, m_r$
  be signed integers, and set $D \assign \{ (\zeta_1, m_1), \ldots, (\zeta_r,
  m_r) \}$. Compute the set $\mathcal{L} (D)$ of rational functions $G / H \in
  \mathbb{K} (\mathcal{C})$ such that the valuation of $G / H$ at $\zeta_i$ is
  $\geqslant - m_i$ for $i = 1, \ldots, r$, and that for any other point
  $\zeta$ on $\mathcal{C}$ the valuation of $G / H$ at $\zeta$ is $\geqslant
  0$.
\end{descriptioncompact}

Problem~B corresponds to the particular case of Problem~C where $F (x, y, z) =
y$, so $\mathcal{C}$ is a projective line in the projective plane. However
Problem~C turns out to be much more difficult in general, and requires several
algebraic constructions. The set $D$ represents a {\tmem{divisor}}
of~$\mathcal{C}$; see definition in Section~\ref{s:divisor}. $\mathcal{L} (D)$
is the {\tmem{Riemann--Roch space}} of $D$.

In their 1874 paper~{\cite{BrillNoether1874}} Brill and Noether designed a
method to compute a $\mathbb{K}$-basis of $\mathcal{L} (D)$, that extends the
above elementary ideas: one first computes a common denominator $H \in
\mathbb{K} [x, y, z]$ of $\mathcal{L} (D)$ and then one deduces a numerator
basis via linear algebra. The computation of $H$ depends on the singular locus
of $\mathcal{C}$. Presenting the Brill and Noether method for computing
Riemann--Roch spaces is the central goal of this paper; see
Section~\ref{s:BNmethod}.

\subsection{Application to error correction}

For applications it is usual to take $\mathbb{K}=\mathbb{F}_q$, the finite
field with $q$ elements. Reed--Solomon error-correcting codes are a popular
technique to represent data in the form of vectors such that the data can be
recovered even if some vector coordinates are corrupted during transmission or
storage of the data. Let $k < n$ be integers, let $\alpha_1, \ldots, \alpha_n$
be distinct points in~$\mathbb{K}$. Let $(a_1, \ldots, a_k) \in \mathbb{K}^k$
be the data to encode. We first interpolate the unique polynomial $f \in
\mathbb{K} [x]_{< k}$ such that $f (\alpha_i) = a_i$ for $i = 1, \ldots, k$,
then we encode the data as
\[ (f (\alpha_1), \ldots, f (\alpha_n)) . \]
This {\tmem{encoding}} can be regarded as a redundant representation of $(a_1,
\ldots, a_k)$: even if a few values $f (\alpha_i)$ are lost, one can recover
$f$ and therefore $(a_1, \ldots, a_k)$. This recovering task is called
{\tmem{decoding}}. The {\tmem{defining vector space}} of this Reed--Solomon
code is the image of the injective map
\begin{eqnarray*}
  \mathcal{L} \left( \left\{ \left( (1 : 0), \, k - 1 \right) \right\} \right)
  & \longhookrightarrow & \mathbb{K}^n\\
  h (x, z) & \longmapsto & (h (\alpha_1, 1), \ldots, h (\alpha_n, 1)) .
\end{eqnarray*}
As seen in Example~\ref{ex:rr-rs}, functions of $\mathcal{L} \left( \left\{
\left( (1 : 0), \, k - 1 \right) \right\} \right)$ are in one-to-one
correspondence with $\mathbb{K} [x]_{< k}$. Unfortunately they suffer from a
limitation: the cardinality of~$\mathbb{K}$ must be~$\geqslant n$.

So-called {\tmem{Algebraic Geometry}} (AG) codes are a generalization of
Reed--Solomon codes which enjoys the same properties with less limitation.
They are based on the Riemann--Roch spaces introduced in Problem~C. Let $G_1 /
H, \ldots, G_{\ell} / H$ represent a basis of $\mathcal{L} (D)$, let $\xi_1,
\ldots, \xi_n$ be distinct points on $\mathcal{C}$ disjoint from $\zeta_1,
\ldots, \zeta_r$, and let
\begin{eqnarray*}
  \mathcal{E}_n : \quad \mathcal{L} (D) & \longrightarrow & \mathbb{K}^n\\
  F / G & \longmapsto & ((F / G) (\xi_1), \ldots, (F / G) (\xi_n)) .
\end{eqnarray*}
With $\bar{D} \assign D \cup \{ (\xi_1, - 1), \ldots, (\xi_n, - 1) \}$, the
kernel of $\mathcal{E}_n$ is $\mathcal{L}_{\mathcal{C}} (\bar{D})$. The
dimensions of $\mathcal{L} (D)$ and $\mathcal{L} (\bar{D})$ can be estimated
using the Riemann--Roch theorem~{\cite[Chapter~1,
Section~5]{Stichtenoth2009}}. We will not enter into these details here. The
favorable situation for correcting codes occurs when there exists an integer
$k < n$ such that the map $\mathcal{E}_k$ is one-to-one. If $\mathbb{K}$ is a
finite field $\mathbb{F}_q$, then the maximum number $n$ of evaluation points
is the number of points of $\mathcal{C}$ with coordinates in $\mathbb{F}_q$,
called the {\tmem{rational points}} of the curve. An algebraic curve can have
up to $q + 1 + 2 g \sqrt{q}$ rational points over the finite field
$\mathbb{F}_q$, where $g$ is a nonnegative number called the \tmtextit{genus}
of $\mathcal{C}$; see also Remark~\ref{rk:genus}. There exist curves
$\mathcal{C}$, said to be \tmtextit{maximal}, whose number of rational points
attains this upper bound. Such curves, but also any curve with many rational
points, allow to construct AG codes longer than the Reed--Solomon ones over
the same finite field. AG codes became particularly famous after the work of
Tsfasman, Vl{\u a}du{\c t}, and Zink~{\cite{TsfasmanVladutZink1982}} proving
that some sequences of codes built from modular curves have better asymptotic
performances than random codes.

The construction of AG codes is one of the major motivations for computing
Riemann--Roch spaces. Nowadays AG codes are used in information theory
protocols, such as secure multi-party computation and zero-knowledge proofs;
see~{\cite{BordageLhotelNardiRandriam2022,CramerRambaudXing2021}}, they can be
used to assert resilience in distributed storage systems (\tmtextit{e.g.~clouds})~{\cite{Barg2017}} or to build families of \tmtextit{quantum
error-correcting codes} with good parameters~{\cite{quantum-AG}}.
See~{\cite{Couvreur2021}} for a recent survey on AG codes and their
applications. Of course, Riemann--Roch spaces are useful in classical contexts
too. For instance, these spaces occur in arithmetic operations on Jacobian of
curves~{\cite{KhuriMakdisi2007}}, in number theory, and algebraic geometry. To
end this section we briefly explain how Riemann--Roch spaces can be used to
share secrets.

\subsection{Application to secret sharing}

A {\tmem{secret sharing scheme}} is a way for $n$ players to carry parts of a
secret message $\mathcal{S}$ such that the knowledge of all the parts allows
to reconstruct the original message, but no proper subset of these parts leaks
anything about $\mathcal{S}$. A secret sharing scheme with threshold $t$ is a
variation of the latter where the knowledge of $t$ parts of the message is
enough to obtain it, but any subset of $t - 1$ parts discloses no information
at all on $\mathcal{S}$.

Shamir's secret sharing scheme~{\cite{Shamir1979}} provides an elegant
construction based on Lagrange interpolation. The message $\mathcal{S}$ is
given as an element of a finite field $\mathbb{F}_q$, say~$a_0$. Then, one
draws a random polynomial $f$ in $\mathbb{F}_q [x]_{< t}$ with constant
coefficient $a_0$, and distinct points $\alpha_1, \ldots, \alpha_n$ in
$\mathbb{F}_q$. The shared parts of the secret are $(\alpha_1, f (\alpha_1)),
\ldots, (\alpha_n, f (\alpha_n))$. Hence, any subset of $t$ players can
recover the message by interpolation.

Shamir's scheme shares multiple properties with Reed--Solomon codes, including
the limitation on the number of players. Replacing polynomials with rational
functions in a Riemann--Roch space, and elements over a finite field with
rational points on an algebraic curve, allows to reduce drastically the size
of the field with respect to the number of players~$n$. This idea was
introduced in~{\cite{ChenCramer2006}} and then developed
in~{\cite{CascudoChenCramerXing2009,CascudoCramerXing2011,CramerRambaudXing2021}}.
For instance, using maximal elliptic curves (that is when the genus of the
curve is $g = 1$) over a finite field $\mathbb{F}_q$ allows to gain $2
\sqrt{q} - 1$ players more with respect to Shamir's classical scheme.

\begin{remark}
  Actually, when replacing Reed--Solomon codes by AG ones, the secret sharing
  thresholds behave slightly differently. Precisely, an AG code in
  $\mathbb{F}_q^n$ yields a secret sharing scheme with two thresholds $t_1
  \geqslant t_2$ such that any coalition of $t_1$ players can recover the
  secret while any coalition of $t_2$ players or less cannot get any
  information on the secret. In the Reed--Solomon setting the two thresholds
  are consecutive, while they are not in the AG setting;
  see~{\cite[Lemma~8]{couvreur:hal-03400779}} for further details.
\end{remark}

\section{Prerequisites}\label{s:preliminaries}

This section gathers usual definitions and constructions related to zero and
one dimensional Zariski closed sets, in order to introduce plane algebraic
curves and study their intersection. The proofs are not all repeated. Instead,
we refer the reader to textbooks. In what follows, we let $\mathbb{K}$ be a
field and $\bar{\mathbb{K}}$ be a fixed algebraic closure of it.

\subsection{Zariski closed sets}

The affine space of dimension $n $ over $\bar{\mathbb{K}}$ is written
$\mathbb{A}^n$. For a set $S$ of polynomials in $\mathbb{K} [x_1, \ldots,
x_n]$, we write $\mathcal{V}_{\mathbb{A}} (S)$ for the {\tmem{Zariski closed
set}} in $\mathbb{A}^n$ defined as the common zeros of the elements of $S$,
that is
\[ \mathcal{V}_{\mathbb{A}} (S) \assign \{ \zeta \in \mathbb{A}^n \suchthat f
   (\zeta) = 0, \forall f \in S \} . \]
If $\mathcal{V}_{\mathbb{A}} (S)$ is finite and non-empty, then it is said to
have {\tmem{dimension zero}}. An ideal $I$ is said to be {\tmem{zero
dimensional}} if the dimension of $\mathcal{V}_{\mathbb{A}} (I)$ is zero.

The projective space of dimension $n $ over $\bar{\mathbb{K}}$ is denoted by
$\mathbb{P}^n$. For a set $S$ of homogeneous polynomials in $\mathbb{K} [x_0,
\ldots, x_n]$, we write $\mathcal{V}_{\mathbb{P}} (S)$ for the {\tmem{Zariski
closed set}} in $\mathbb{P}^n$ defined as the common zeros of the elements of
$S$, that is
\[ \mathcal{V}_{\mathbb{P}} (S) \assign \{ \zeta \in \mathbb{P}^n \suchthat F
   (\zeta) = 0, \forall F \in S \} . \]
If $\mathcal{V}_{\mathbb{P}} (S)$ is finite and non-empty, then it is said to
have {\tmem{dimension zero}}. For a non-constant homogeneous polynomial $F \in
\mathbb{K} [x, y, z]$ the set $\mathcal{V}_{\mathbb{P}} (F)$ is called an
{\tmem{algebraic curve}}. Curves are {\tmem{one dimensional}} Zariski closed
sets.

If $\mathbb{M} \assign \mathbb{K} [x_1, \ldots, x_n]$ is a polynomial ring and
$\zeta$ a point in $\mathbb{K}^n$, then $\mathbb{M}_{\zeta}$ represents the
{\tmem{local ring}} of the rational functions $A / B$ in $\mathbb{K} (x_1,
\ldots, x_n)$ such that $B (\zeta) \neq 0$. The following classical result
will be needed for $n = 2$ variables.

\begin{proposition}
  \label{pp:zero-dim-ideal}{\tmem{({\cite[Chapter~2, Section~9,
  Proposition~6]{FultonBook}}, {\cite[Chapter~4,
  Section~2]{CoxLittleOShea2005}}, or
  {\cite[Chapter~4]{KreuzerRobbiano2016}})}} Assume that $\mathbb{K}$ is
  algebraically closed. Let $I$ be a zero dimensional ideal in $\mathbb{M}
  \assign \mathbb{K} [x_1, \ldots, x_n]$. Then, we have
  \[ \mathbb{M}/ I \cong \bigoplus_{\zeta \in \mathcal{V}_{\mathbb{A}} (I)}
     \mathbb{M}_{\zeta} / (I\mathbb{M}_{\zeta}), \]
  where each summand is a local $\mathbb{K}$-algebra of finite dimension.
\end{proposition}

\subsection{Resultant}\label{s:resultant}

Let $\mathbb{D}$ be an integral domain. The {\tmem{resultant}} of two
polynomials $f$ and $g$ in $\mathbb{D} [x]$ of respective degrees~$m$ and~$n$
is defined as the determinant of the map
\begin{eqnarray*}
  \mathbb{D} [x]_{< n} \times \mathbb{D} [x]_{< m} & \longrightarrow &
  \mathbb{D} [x]_{< m + n}\\
  (u, v) & \longmapsto & uf + vg,
\end{eqnarray*}
and is written $\tmop{Res} (f, g)$. By construction, $\tmop{Res} (f, g)$
belongs to the ideal $(f, g)$. If $f, g$ are multivariate polynomials then
their resultant with respect to an indeterminate~$\mathit{x}$ is denoted by
$\tmop{Res}_x (f, g)$. We will need the following propositions.

\begin{proposition}
  \label{pp:res-premier}{\tmem{({\cite[Chapitre~6,
  Corollaire~6.3]{aecf-2017-livre}} or {\cite[Chapter~3, Section~6,
  Proposition~3]{CoxLittleOShea2013}})}} Two polynomials $f, g$ in $\mathbb{D}
  [x]$ have no common factor of degree $\geqslant 1$ if, and only if,
  $\tmop{Res} (f, g) \neq 0$.
\end{proposition}

\begin{proposition}
  \label{pp:res-charpoly}{\tmem{({\cite[Chapitre~6,
  Lemme~6.9]{aecf-2017-livre}} or {\cite[Chapter~3,
  Proposition~1.5]{CoxLittleOShea2005}}}}) If $f$ is monic, then $\tmop{Res}
  (f, g)$ is the determinant of the multiplication endomorphism by $g$ in
  $\mathbb{K} [x] / (f (x))$.
\end{proposition}

\begin{proposition}
  \label{pp:res-mul}{\tmem{{\cite[Chapitre~6,
  Corollaire~6.6]{aecf-2017-livre}}}} For all polynomials $f$, g, and $h$ in
  $\mathbb{D} [x]$ we have
  \[ \tmop{Res} (fg, h) = \tmop{Res} (f, h) \tmop{Res} (g, h) . \]
\end{proposition}

{\noindent}The resultant is a long-established tool for polynomial system
solving. Let $F, G \in \mathbb{K} [x, y, z]$ be two non-constant homogeneous
coprime polynomials, then the intersection of the curves defined by $F$ and
$G$ is $\mathcal{V}_{\mathbb{P}} (F) \cap \mathcal{V}_{\mathbb{P}} (G)
=\mathcal{V}_{\mathbb{P}} (F, G)$. We now recall a simple method for computing
curve intersection, that makes use of the resultant of $F$ and $G$ regarded in
$\mathbb{K} [x, z] [y]$. We set
\[ R (x, z) \assign \tmop{Res}_y (F (x, y, z), G (x, y, z)), \]
and introduce specific technical assumptions:

\begin{descriptioncompact}
  \item[A\tmrsub{1}] $\deg_y F = \deg F$ and $F$ is monic in $y$;
  {\nopagebreak}
  
  \item[A\tmrsub{2}] $\deg_x R (x, z) = \deg (R (x, z))$,
\end{descriptioncompact}

\noindent where $\deg_y F$ and $\deg_x R$ represent the partial degrees of $F$ in $y$
and of $R$ in $x$.

\begin{lemma}
  \label{lm:generic-intersection}Let $F$ and $G$ be non-constant coprime
  homogeneous polynomials in $\mathbb{K} [x, y, z]$ such that A\tmrsub{1} and
  A\tmrsub{2} hold. Then, the following assertions hold:
  \begin{enumerateroman}
    \item \label{it:curve-inter-a}$\mathcal{V}_{\mathbb{P}} (F, G)$ is
    contained in the affine chart defined by $z = 1$;
    
    \item \label{it:curve-inter-b}$\mathcal{V}_{\mathbb{A}} (F (x, y, 1), G
    (x, y, 1)) = \bigcup_{R (\zeta_x, 1) = 0} \{ (\zeta_x, \zeta_y) \suchthat
    F (\zeta_x, \zeta_y, 1) = G (\zeta_x, \zeta_y, 1) = 0 \}$;
    
    \item \label{it:curve-inter-c}$\mathcal{V}_{\mathbb{A}} (F (x, y, 1), G
    (x, y, 1))$ is finite and non-empty.
  \end{enumerateroman}
\end{lemma}

\begin{proof}
  Let $(\zeta_x : \zeta_z) \in \mathbb{P}^1$, and let $M (x, z)$ denote the
  matrix of the multiplication endomorphism by $G (x, y, z)$ in $\mathbb{K}
  (x, z) [y] / (F (x, y, z))$ in the basis $1, y, \ldots, y^{\deg_y F - 1}$.
  Thanks to~A\tmrsub{1}, the entries of $M (x, z)$ belong to $\mathbb{K} [x,
  z]$ and the matrix of the multiplication endomorphism by $G (\zeta_x, y,
  \zeta_z)$ in $\mathbb{K} [y] / (F (\zeta_x, y, \zeta_z))$ in the basis $1,
  y, \ldots, y^{\deg_y F - 1}$ coincides to $M (\zeta_x, \zeta_z)$.
  Proposition~\ref{pp:res-charpoly} implies that
  \[ R (\zeta_x, \zeta_z) = \tmop{Res}_y (F (\zeta_x, y, \zeta_z), G (\zeta_x,
     y, \zeta_z)) . \]
  If $(\zeta_x : \zeta_y : 0) \in \mathcal{V}_{\mathbb{P}} (F, G)$ then by
  Proposition~\ref{pp:res-premier} we have $R (\zeta_x, 0) = 0$,
  hence~A\tmrsub{2} implies $\zeta_x = 0$, so~A\tmrsub{1} yields $\zeta_y =
  0$, whence property~\eqref{it:curve-inter-a}. As
  for~\eqref{it:curve-inter-b}, by Proposition~\ref{pp:res-premier} we know
  that $\zeta_x$ is a root of $R (\zeta_x, 1) = 0$ if, and only if, $F
  (\zeta_x, y, 1)$ and $G (\zeta_x, y, 1)$ share at least one common root
  in~$\bar{\mathbb{K}}$. Finally, thanks to~A\tmrsub{1} the polynomial $F
  (\zeta_x, y, 1)$ is non-zero and hence admits a finite number of roots
  in~$\bar{\mathbb{K}}$, whence~\eqref{it:curve-inter-c}.
\end{proof}

Assumptions~A\tmrsub{1} and~A\tmrsub{2} can be recovered after suitable
changes of coordinates, as stated in the following lemma.

\begin{lemma}
  \label{lm:generic-coords}With the notation as above, let $(\alpha, \beta,
  \gamma) \in \mathbb{K}^3$. If $F (\alpha, 1, \beta) \neq 0$ then the partial
  degree of $F (x + \alpha y, y, z + \beta y)$ in $y$ equals $\deg F$. If $R
  (1, \gamma) \neq 0$ then the partial degree of $R (x, z + \gamma x)$ in~$x$
  equals $\deg R$.
\end{lemma}

\begin{proof}
  Since $F$ is homogeneous, the coefficient of $y^{\deg F}$ in $F (x + \alpha
  y, y, z + \beta y)$ is $F (\alpha, 1, \beta)$. The second assertion about
  $R$ is obtained in the same manner.
\end{proof}

If the cardinality of $\mathbb{K}$ is sufficiently large, then there exists a
triple $(\alpha, \beta, \gamma) \in \mathbb{K}^3$ such that $F (\alpha, 1,
\beta) \neq 0$ and $R (1, \gamma) \neq 0$. This simple fact will be sufficient
in the sequel. Nevertheless, let us mention a stronger assertion, that is a
direct consequence of the classical Schwarz--Zippel
lemma~{\cite[Lemma~6.44]{GathenGerhard2013}}: if $\mathcal{S}$ is a finite
subset of $\mathbb{K}$, then the probability that a triple $(\alpha, \beta,
\gamma)$ taken uniformly at random in $\mathcal{S}^3$ satisfies $F (\alpha, 1,
\beta) R (1, \gamma) \neq 0$ is~$\geqslant 1 - (\deg F + \deg R) / |
\mathcal{S} |$.

\begin{proposition}
  \label{pp:dimension-curve-intersection}Let $F$ and $G$ be non-constant
  homogeneous polynomials in $\mathbb{K} [x, y, z]$. If $F$ and $G$ are
  coprime then $\mathcal{V}_{\mathbb{P}} (F, G)$ is finite and
  non-empty.{\nopagebreak}
\end{proposition}

\begin{proof}
  This is a consequence of Lemmas~\ref{lm:generic-intersection}
  and~\ref{lm:generic-coords} used over $\bar{\mathbb{K}}$ (that is infinite).
\end{proof}

\section{Hensel lemmas}\label{s:hensel}

The usual Hensel lemma concerns polynomials $f \in \mathbb{K} [[x]] [y]$ monic
in $y$: given $g_0$ and $h_0$ in $\mathbb{K} [y]$ monic and coprime such that
$f = g_0 h_0 + O (x)$, then there exist unique monic polynomials $g$ and $h$
in $\mathbb{K} [[x]] [y]$ such that $f = gh$, $g = g_0 + O (x)$ and $h = h_0 +
O (x)$; see~{\cite[Chapter~7, Theorem~7.18]{Eisenbud1995}}
or~{\cite[Chapitre~21, Proposition~21.19]{aecf-2017-livre}} for instance. This
section is devoted to usual extensions of this lemma in the context of
bivariate power series.

\subsection{Bivariate weighted valuations}

Let $\gamma = (\gamma_x, \gamma_y) \in \mathbb{N}^2$ with $\gamma_x \neq 0$,
and let $\tmop{val} (x^a y^b) \assign \gamma_x a + \gamma_y b$ define a
{\tmem{weighted valuation}} over $\mathbb{K} [[x, y]]$. An element in
$\mathbb{K} [[x, y]]$ is said to be {\tmem{quasi-homogeneous}} if its non-zero
terms have the same weighted valuation. It will be convenient to assume that
$\gamma_x$ and $\gamma_y$ are coprime. This makes $\mathbb{K} [[x, y]]$ a
graded ring. The quasi-homogeneous component of
\[ a = \sum_{i, j \geqslant 0} a_{i, j} x^i y^j \in \mathbb{K} [[x, y]] \]
of valuation $e$ (if non-zero) is written
\[ [a]_e \assign \sum_{\gamma_x i + \gamma_y j = e} a_{i, j} x^i y^j \in
   \mathbb{K} [x, y] . \]
The {\tmem{initial form}} of $a$ is denoted by
\[ \tmop{in} (a) \assign [a]_{\tmop{val} (a)}, \]
and we introduce the following notation for truncations of $a$:
\[ [a]_{e ; \eta} \assign \sum_{e \leqslant \gamma_x i + \gamma_y j < e +
   \eta} a_{i, j} x^i y^j \in \mathbb{K} [x, y], \text{ for } \eta > 0. \]
\begin{lemma}
  \label{lm:division}Let $f, g \in \mathbb{K} [[x]] [y]$. If $g$ is monic of
  degree $n$ and if $\tmop{val} (g) = \tmop{val} (y^n)$, then the remainder $f
  \tmop{rem} g$ in the division of $f$ by $g$ satisfies $\tmop{val} (f
  \tmop{rem} g) \geqslant \tmop{val} (f)$.
\end{lemma}

\begin{proof}
  The result is clear whenever $\deg f < n$. Now assume that $m \assign \deg f
  \geqslant n$ and let $f_m (x)$ denote the coefficient of $y^m$ in $f$. We
  let $h \assign f - f_m (x) y^{m - n} g$. By definition we have
  \[ \tmop{val} (f_m (x) y^m) \geqslant \tmop{val} (f) . \]
  Using $\tmop{val} (g) = \tmop{val} (y^n)$ we deduce that $\tmop{val} (f_m
  (x) y^{m - n} g) = \tmop{val} (f_m (x) y^m)$, and then that
  \[ \tmop{val} (h) \geqslant \min (\tmop{val} (f), f_m (x) y^{m - n} g)
     \geqslant \tmop{val} (f) . \]
  Since $\deg h < \deg f$ and $h \tmop{rem} g = f \tmop{rem} g$, the
  conclusion follows from a straightforward induction on $n$.
\end{proof}

\begin{lemma}
  \label{lm:hensel}Let $f \in \mathbb{K} [[x, y]] \setminus \{ 0 \}$ be such
  that $\tmop{in} (f)$ has bounded degree in $y$. Then, there exist unique $u
  \in \mathbb{K} [[x, y]]$ and $g \in \mathbb{K} [[x]] [y]$ such that $f =
  ug$, $g$ is monic with $\deg_y g = \deg_y (\tmop{in} (f))$, and $\tmop{in}
  (u) \in x^m  (\mathbb{K} \setminus \{ 0 \})$ for some $m \in \mathbb{N}$.
\end{lemma}

\begin{proof}
  If $u$ and $g$ exist as specified, then we must have $\tmop{in} (f) =
  \tmop{in} (u) \tmop{in} (g)$. Let us write
  \[ \tmop{in} (f) = \sum_{i \geqslant 0}^n f_{m_i, i} x^{m_i} y^i, \]
  where $\gamma_x m_i + \gamma_y i = \tmop{val} (f)$ for $i = 0, \ldots, n$,
  and with $f_{m_n, n} \neq 0$. Note that $m_i \geqslant m_n$ for {$i = 0,
  \ldots, n$}. Consequently we set $m \assign m_n$, $u_1 \assign f_{m, n} x^m$
  and $g_1 \assign \tmop{in} (f) / u_1$, so $\tmop{in} (u) = u_1$ and
  $\tmop{in} (g) = g_1$ hold necessarily.
  
  For $k \geqslant 1$, assume by induction that there exist $u_k \in
  \mathbb{K} [[x, y]]$ and $g_k \in \mathbb{K} [[x]] [y]$ such that
  \[ [f]_{\tmop{val} (f) ; k} = [u_k g_k]_{\tmop{val} (f) ; k}, \]
  $g_k$ is monic, and $\tmop{in} (g_k) = g_1$. In order to construct $u_{k +
  1}$ and $g_{k + 1}$, we are looking for $\tilde{u} = u_{k + 1} - u_k$ and
  $\tilde{g} = g_{k + 1} - g_k$ quasi-homogeneous of valuation $\tmop{val}
  (u_1) + k$ and $\tmop{val} (g_1) + k$, respectively, and such that
  \[ [f]_{\tmop{val} (f) ; k + 1} = [(u_k + \tilde{u})  (g_k +
     \tilde{g})]_{\tmop{val} (f) ; k + 1}, \]
  and $\deg_y  \tilde{g} < \deg_y g_1$. The latter equation rewrites into
  \[ [f - u_k g_k]_{\tmop{val} (f) ; k + 1} = [u_k  \tilde{g} + \tilde{u}
     g_k]_{\tmop{val} (f) ; k + 1}, \]
  which is further equivalent to
  \[ [f - u_k g_k]_{\tmop{val} (f) + k} = [u_1  \tilde{g} + \tilde{u}
     g_1]_{\tmop{val} (f) + k} . \]
  Since $\deg_y  \tilde{g} < \deg_y g_1$, the polynomial $u_1  \tilde{g}$ must
  be the remainder, written $r$, in the division of $[f - u_k g_k]_{\tmop{val}
  (f) + k}$ by $g_1$, and $\tilde{u}$ is the corresponding quotient. This
  division is well defined because $g_1 = \tmop{in} (f) / u_1$ is monic in
  $y$. Lemma~\ref{lm:division} implies $\tmop{val} (r) \geqslant \tmop{val}
  (f) + k$. Consequently, $\tmop{val} (r) > \tmop{val} (f) = \tmop{val} (x^m
  y^n)$. The inequality $\deg_y r < n$ implies that $\tmop{val}_x (r) > m$, so
  $r$ is a multiple of $u_1$, and therefore $\tilde{g}$ does exist and is
  given by $\tilde{g} \assign r / u_1$. In order to conclude the proof it
  suffices to take $u \assign \lim_{k \rightarrow \infty} u_k$ and $g \assign
  \lim_{k \rightarrow \infty} g_k$.
\end{proof}

The following lemma is sometimes called the {\tmem{Weierstra{\ss}
normalization}}, or {\tmem{Weierstra{\ss} preparation theorem}} in the
analytic context.

\begin{lemma}
  \label{lm:weierstrass}Let $f = \sum_{i = 0}^{\infty} f_i (x) y^i \in
  \mathbb{K} [[x]] [[y]]$ be such that there exists a positive integer $n$
  satisfying $\tmop{val}_x (f_n (x)) = 0$ and $\tmop{val}_x (f_i (x)) > 0$ for
  $i < n$. Then, there exist $u$ invertible in $\mathbb{K} [[x, y]]$ and $g$
  monic of degree $n$ in $\mathbb{K} [[x]] [y]$ such that $f = ug$.
\end{lemma}

\begin{proof}
  With the notation as above, set $\gamma_x \assign n$ and $\gamma_y \assign
  1$. Clearly $\tmop{val} (f) = n$ and $\tmop{in} (f)$ has degree $n$ in $y$.
  We have $\tmop{in} (g) = \tmop{in} (f) / f_n (0)$ and the conclusion follows
  from Lemma~\ref{lm:hensel}.
\end{proof}

\subsection{Hensel lifting for polynomials}

In the preceding lemmas we have seen how factorization in $\mathbb{K} [[x,
y]]$ reduces to factorization in $\mathbb{K} [[x]] [y]$. The next lemma
addresses the latter case for weighted valuations.

\begin{lemma}
  \label{lm:hensel-bis}Let $\tmop{val}$ denote a weighted valuation as above,
  and let $f \in \mathbb{K} [[x]] [y]$ be monic of degree~$n$ such that
  $\tmop{val} (f) = \tmop{val} (y^n)$. Assume that $\tmop{in} (f)$ factorizes
  into two quasi-homogeneous coprime monic polynomials $g_1$ and $h_1$. Then,
  there exist monic polynomials $g$ and $h$ in $\mathbb{K} [[x]] [y]$ such
  that $f = gh$, $\tmop{in} (g) = g_1$, and $\tmop{in} (h) = h_1$.
\end{lemma}

\begin{proof}
  The B{\'e}zout relation for $g_1$ and $h_1$ regarded in $\mathbb{K} ((x))
  [y]$ can be written
  \begin{equation}
    ug_1 + vh_1 = x^m, \label{eq:hensel-bis-i}
  \end{equation}
  where $m \in \mathbb{N}$, $u$ and $v$ belong to $\mathbb{K} [[x]] [y]$ and
  satisfy $\deg_y u < \deg_y h_1$, $\deg_y v < \deg_y g_1$. Since $g_1$ and
  $h_1$ are quasi-homogeneous we may even take $u$ and $v$ quasi-homogeneous
  in~\eqref{eq:hensel-bis-i}, so we have
  \begin{equation}
    \tmop{val} (u) + \tmop{val} (g_1) = \tmop{val} (v) + \tmop{val} (h_1) =
    \gamma_x m. \label{eq:hensel-bis-ii}
  \end{equation}
  For $k \geqslant 1$, by induction we may assume that there exist $g_k$ and
  $h_k$ monic in $\mathbb{K} [[x]] [y]$ such that
  \[ [f]_{\tmop{val} (f) ; k} = [g_k h_k]_{\tmop{val} (f) ; k}, \]
  and $\tmop{in} (g_k) = g_1$. In order to construct $g_{k + 1}$ and $h_{k +
  1}$, we are looking for $\tilde{g} = g_{k + 1} - g_k$ and $\tilde{h} = h_{k
  + 1} - h_k$ quasi-homogeneous of valuation $\tmop{val} (g_1) + k$ and
  $\tmop{val} (h_1) + k$, respectively, such that
  \[ [f]_{\tmop{val} (f) ; k + 1} = [(g_k + \tilde{g})  (h_k +
     \tilde{h})]_{\tmop{val} (f) ; k + 1}, \]
  $\deg_y  \tilde{g} < \deg_y g_1$, and $\deg_y  \tilde{h} < \deg_y h_1$. The
  latter equation rewrites into
  \[ [f - g_k h_k]_{\tmop{val} (f) ; k + 1} = [g_k  \tilde{h} + \tilde{g}
     h_k]_{\tmop{val} (f) ; k + 1}, \]
  which is equivalent to
  \[ [f - g_k h_k]_{\tmop{val} (f) + k} = g_1  \tilde{h} + \tilde{g} h_1 . \]
  Thanks to the usual Chinese remaindering theorem over $\mathbb{K} ((x))$,
  and using~\eqref{eq:hensel-bis-i}, $\tilde{g}$ and $\tilde{h}$ must satisfy
  \begin{eqnarray*}
    \tilde{g} & = & ((v / x^m)  [f - g_k h_k]_{\tmop{val} (f) + k} \tmop{rem}
    g_1),\\
    \tilde{h} & = & ((u / x^m)  [f - g_k h_k]_{\tmop{val} (f) + k} \tmop{rem}
    h_1),
  \end{eqnarray*}
  where $p \tmop{rem} q$ represents the remainder in the division of $p$ by
  $q$. Lemma~\ref{lm:division} implies
  \[ \tmop{val} (v [f - g_k h_k]_{\tmop{val} (f) + k} \tmop{rem} g_1)
     \geqslant \tmop{val} (v [f - g_k h_k]_{\tmop{val} (f) + k}) > \tmop{val}
     (v) + \tmop{val} (f) . \]
  Combined with~\eqref{eq:hensel-bis-ii} we deduce that
  \begin{equation}
    \tmop{val} (v [f - g_k h_k]_{\tmop{val} (f) + k} \tmop{rem} g_1) >
    \tmop{val} (g_1) + \gamma_x m. \label{eq:hensel-bis-iii}
  \end{equation}
  Then, thanks to
  \[ \deg_y (v [f - g_k h_k]_{\tmop{val} (f) + k} \tmop{rem} g_1) < \deg_y
     g_1, \]
  a monomial $x^a y^b$ with non-zero coefficient in $v [f - g_k
  h_k]_{\tmop{val} (f) + k} \tmop{rem} g_1$ satisfies
  \begin{equation}
    \gamma_y b \leqslant \gamma_y \deg_y g_1 = \tmop{val} (g_1) .
    \label{eq:hensel-bis-iv}
  \end{equation}
  Combining~\eqref{eq:hensel-bis-iii} and~\eqref{eq:hensel-bis-iv} yields
  \[ \gamma_x a + \tmop{val} (g_1) \geqslant \gamma_x a + \gamma_y b >
     \tmop{val} (g_1) + \gamma_x m. \]
  Since $\gamma_x \neq 0$, it follows that $a \geqslant m$ and that the above
  value found for $\tilde{g}$ actually belongs to $\mathbb{K} [[x]] [y]$. In
  the same way we obtain that $\tilde{h} \in \mathbb{K} [[x]] [y]$. In order
  to conclude the proof it suffices to take $g \assign \lim_{k \rightarrow
  \infty} g_k$ and $h \assign \lim_{k \rightarrow \infty} h_k$.
\end{proof}

\subsection{Newton polygon and irreducibility}

\begin{figure}[t]
  \begin{center}
    \includegraphics[width=9.3255772005772cm,height=4.95822674799948cm]{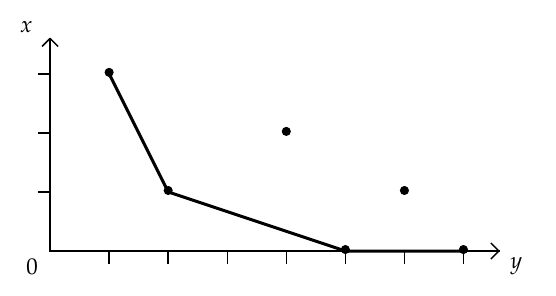}
  \end{center}
  \caption{\label{fig:newton-poly}Newton polygon of $f = x^3 y + 2 xy^2 - x^2
  y^4 + y^5 + 3 xy^6 + y^7 \in \mathbb{Q} [[x]] [y]$.}
\end{figure}

The Newton polygon is a tool to study the behavior of roots of polynomials.
Let
\[ f = \sum_{i = 0}^n f_i (x) y^i \in \mathbb{K} ((x)) [y] \]
be of degree $n$. The {\tmem{Newton polygon}} of $f$ is the {\tmem{lower
border of the convex hull}} in $\mathbb{R}_{\geqslant 0} \times \mathbb{R}$ of
the set of points
\[ \{ (i, \tmop{val}_x (f_i)) \suchthat 0 \leqslant i \leqslant n, f_i \neq 0
   \} . \]
Figure~\ref{fig:newton-poly} illustrates this definition. The Newton polygon is a broken line with edges $(i_0, j_0),
\ldots, (i_r, j_r)$, with $r \geqslant 0$, such that:
\begin{enumerate}
  \item $i_0 < i_1 < \cdots < i_r$,
  
  \item $j_k = \tmop{val} (f_{i_k})$ for $k = 0, \ldots, r$,
  
  \item \label{enum:newton-poly}for all $i = 0, \ldots, n$ such that $f_i \neq
  0$ the point $(i, \tmop{val}_x (f_i))$ is not strictly below the Newton
  polygon.
\end{enumerate}
If $f$ is zero then its Newton polygon is empty. If $f$ has a single term then
its Newton polygon reduces to a single vertex, and it is said to be
{\tmem{degenerate}}. Otherwise, any {$k = 1, \ldots, r$} determines an edge
$E_k$ with vertices $(i_{k - 1}, j_{k - 1})$ and $(i_k, j_k)$, and of
{\tmem{slope}} $(j_k - j_{k - 1}) / (i_k - i_{k - 1})$. The {\tmem{Newton
polynomial}} associated to $E_k$ is
\[ \sum_{i \in \mathbb{N}, (i, j) \in E_k} [f_i (x)]_j y^{i - i_{k - 1}} \in
   \mathbb{K} (x) [y], \]
where $[f_i (x)]_j$ represents the term of degree $j$ in $f_i$.

\begin{example}
  The edge $E_2$ of vertices $(2, 1)$ and $(5, 0)$ in
  Figure~\ref{fig:newton-poly} has associated Newton polynomial $2 x + y^3$.
\end{example}

Let $f = \sum_{i = 0}^n f_i (x) y^i \in \mathbb{K} [[x]] [y]$ be monic of
degree $n \geqslant 1$, and let $(i_0, j_0), \ldots, (i_r, j_r)$ represent the
Newton polygon of $f$. If this polygon is degenerate then $f$ is reducible
unless $n = 1$. We draw a useful consequence of Lemma~\ref{lm:hensel}.

\begin{lemma}
  \label{lm:irreducible}Let $f = \sum_{i = 0}^n f_i (x) y^i \in \mathbb{K}
  ((x)) [y]$ be monic of degree $n \geqslant 1$, irreducible, and such that
  $f_0 (x) \neq 0$ and $\tmop{val}_x (f_0) \geqslant 0$. Then, the Newton
  polygon of $f$ has a single edge of vertices $(0, \tmop{val}_x (f_0))$ and
  $(n, 0)$.
\end{lemma}

\begin{proof}
  Let $k \geqslant 0$ be the smallest integer such that $g \assign x^k f$
  belongs to $\mathbb{K} [[x]] [y]$, and let $(0, \tmop{val}_x (x^k f_0))$ and
  $(l, \tmop{val}_x (x^k f_l))$ represent the vertices of the first edge $E$
  of the Newton polygon of $g$. There exist coprime integers $\gamma_x
  \geqslant 1$, $\gamma_y \geqslant 0$ such that
  \[ \frac{\gamma_y}{\gamma_x} = \frac{\tmop{val}_x (f_0) - \tmop{val}_x
     (f_l)}{l} . \]
  If we had $l < n$, then the Newton polygon of $g$ would have at least two
  edges, so Lemma~\ref{lm:hensel} (used for $g$) would yield a non-trivial
  factor $h$ of $g$ of degree $l$ in $y$. Since $f$ is irreducible this is not
  possible, whence $l = n$, $k = 0$, and $\tmop{in} (f)$ is the Newton
  polynomial of $E$.
\end{proof}

\begin{proposition}
  \label{pp:irreducible}Let $f = \sum_{i = 0}^n f_i (x) y^i \in \mathbb{K}
  [[x]] [y]$ be monic of degree $n \geqslant 1$, irreducible, and such that
  $f_0 (x) \neq 0$ and $f_0 (0) = 0$. Then, the Newton polygon of $f$ has a
  single edge whose Newton polynomial writes
  \[ x^{\tmop{val}_x (f_0)} \theta \left( \frac{y^{\gamma_x}}{x^{\gamma_y}}
     \right) \]
  where $\theta \in \mathbb{K} [t]$ is monic and irreducible of degree $d = n
  / \gamma_y$, and $\gamma_x \geqslant 1$, $\gamma_y \geqslant 1$ are coprime
  integers such that
  \[ \frac{\gamma_y}{\gamma_x} = \frac{\tmop{val}_x (f_0)}{n} . \]
\end{proposition}

\begin{proof}
  From Lemma~\ref{lm:irreducible} we know that the Newton polygon of $f$ has a
  single edge $E$ and that $\tmop{in} (f)$ is the Newton polynomial of $E$.
  For a monomial $x^a y^b$ of weighted-valuation
  \[ \tmop{val} (f) = \gamma_x a + \gamma_y b = \gamma_y n = \gamma_x
     \tmop{val}_x (f_0) \]
  there exists a unique integer $k \geqslant 0$ such that $\tmop{val}_x (f_0)
  - a = k \gamma_y$ and $b = k \gamma_x$. We obtain:
  \[ x^a y^b = x^{\tmop{val}_x (f_0) - k \gamma_y} y^{k \gamma_x} =
     x^{\tmop{val}_x (f_0)}  \left( \frac{y^{\gamma_x}}{x^{\gamma_y}}
     \right)^k, \]
  hence $\tmop{in} (f)$ writes
  \[ \tmop{in} (f) = x^{\tmop{val}_x (f_0)} \theta \left(
     \frac{y^{\gamma_x}}{x^{\gamma_y}} \right), \]
  where $\theta$ is a monic polynomial in $\mathbb{K} [t]$. If $\theta$ had a
  non-trivial irreducible factor $\tilde{\theta}$ then
  \[ x^{\gamma_y \deg \tilde{\theta}}  \tilde{\theta} \left(
     \frac{y^{\gamma_x}}{x^{\gamma_y}} \right) \]
  would be a non-trivial factor of $\tmop{in} (f)$. This is not possible
  because of Lemma~\ref{lm:hensel-bis}, so $\theta$ is irreducible.
\end{proof}

\section{Valuations}\label{s:valuations}

The notion of \tmtextit{valuation} has been introduced just before
Section~\ref{s:intro-rr}. The goal of this section is to highlight how
valuations extend in algebraic extensions, and to introduce the notion of
\tmtextit{uniformizing parameter}. If $f$ is irreducible in $\mathbb{K} [[x,
y]]$, then the valuations on $\mathbb{K} [[x, y]] / (f)$ are in one-to-one
correspondence to those on $\tmop{Frac} (\mathbb{K} [[x, y]] / (f))$.

\subsection{Algebraic extensions}

The valuation $\tmop{val}_x$ is the unique valuation on $\mathbb{K} ((x))$
normalized such that $\tmop{val}_x (x) = 1$. We examine how this valuation
extends in algebraic extensions.

\begin{lemma}
  \label{lm:slope}Let $f \in \mathbb{K} [[x]] [y]$ be irreducible, monic of
  degree $n$ in $y$, and such that $f (0, 0) = 0$. If~$v$ is a valuation on
  $\mathbb{K} [[x, y]] / (f)$ that extends $\tmop{val}_x$, then we have
  \[ v (y \tmop{mod} f) = \frac{\tmop{val}_x (f (x, 0))}{n}, \]
  where $y \tmop{mod} f$ denotes the class of $y$ in $\mathbb{K} [[x, y]] /
  (f)$.
\end{lemma}

\begin{proof}
  If $f = y$ then the statement clearly holds, so let us assume that $f \neq
  y$. Let us write $f = \sum_{i = 0}^n f_i (x) y^i$. Since $f$ is monic in $y$
  we have $\tmop{val}_x (f_n) = 0$. Since $f$ is irreducible we have $f (x, 0)
  = f_0 (x) \neq 0$. By Proposition~\ref{pp:irreducible}, the Newton polygon
  of $f$ has a single edge whose vertices are $(0, \tmop{val}_x (f_0))$ and
  $(n, 0)$. This entails that
  \[ \tmop{val}_x (f_i) \geqslant \frac{\tmop{val}_x (f_0)}{n}  (n - i) \]
  for $i = 0, \ldots, n$; see property~\eqref{enum:newton-poly} of the
  definition of the Newton polygon. If a valuation $v$ on $\mathbb{K} [[x, y]]
  / (f)$ extends $\tmop{val}_x$, then we necessarily have
  \begin{eqnarray*}
    nv (y \tmop{mod} f) & = & v ((y^n - f) \tmop{mod} f)\\
    & = & v \left( \sum_{i = 0}^{n - 1} f_i (x) y^i \tmop{mod} f \right)\\
    & \geqslant & \min_{i = 0, \ldots, n - 1} (v (f_i (x) y^i \tmop{mod}
    f))\\
    & = & \min_{i = 0, \ldots, n - 1} (\tmop{val}_x (f_i) + iv (y \tmop{mod}
    f))\\
    & \geqslant & \min_{i = 0, \ldots, n - 1} \left( \frac{\tmop{val}_x
    (f_0)}{n}  (n - i) + iv (y \tmop{mod} f) \right) .
  \end{eqnarray*}
  Let $j$ be a value of $i \in \{ 0, \ldots, n - 1 \}$ for which
  $\frac{\tmop{val}_x (f_0)}{n}  (n - i) + iv (y \tmop{mod} f)$ is minimal, so
  we have
  \[ nv (y \tmop{mod} f) \geqslant \frac{\tmop{val}_x (f_0)}{n}  (n - j) + jv
     (y \tmop{mod} f), \]
  that implies
  \[ (n - j) v (y \tmop{mod} f) \geqslant \frac{\tmop{val}_x (f_0)}{n}  (n -
     j), \]
  whence
  \begin{equation}
    v (y \tmop{mod} f) \geqslant \frac{\tmop{val}_x (f_0)}{n} .
    \label{eq:slope1}
  \end{equation}
  On the other hand, we verify that
  \begin{eqnarray*}
    \tmop{val}_x (f_0) & = & v ((f - f_0 (x)) \tmop{mod} f)\\
    & = & v \left( \sum_{i = 1}^n f_i (x) y^i \tmop{mod} f \right)\\
    & \geqslant & \min_{i = 1, \ldots, n} (v (f_i (x) y^i \tmop{mod} f))\\
    & = & \min_{i = 1, \ldots, n} (\tmop{val}_x (f_i) + iv (y \tmop{mod} f))
    .\\
    & \geqslant & \min_{i = 1, \ldots, n} \left( \frac{\tmop{val}_x (f_0)}{n}
    (n - i) + iv (y \tmop{mod} f) \right) .
  \end{eqnarray*}
  Again, let $j$ be a value of $i \in \{ 1, \ldots, n \}$ for which
  $\frac{\tmop{val}_x (f_0)}{n}  (n - i) + iv (y \tmop{mod} f)$ is minimal.
  Then we obtain
  \[ \tmop{val}_x (f_0) \geqslant \frac{\tmop{val}_x (f_0)}{n}  (n - j) + jv
     (y \tmop{mod} f), \]
  whence
  \begin{equation}
    v (y \tmop{mod} f) \leqslant \frac{\tmop{val}_x (f_0)}{n} .
    \label{eq:slope2}
  \end{equation}
  The combination of~\eqref{eq:slope1} and~\eqref{eq:slope2} concludes the
  proof.
\end{proof}

Let $f \in \mathbb{K} [[x]] [y]$ be a monic irreducible polynomial of degree
$n$. Let $a \in \mathbb{K} ((x)) [y] / (f)$ be represented by $A \in
\mathbb{K} ((x)) [y]_{< n}$. We regard $\mathbb{K} ((x)) [y] / (f)$ as a
$\mathbb{K} ((x))$-algebra of dimension $n$. The {\tmem{minimal polynomial}}
of $a$ is the monic polynomial $\mu \in \mathbb{K} ((x)) [t]$ of smallest
degree in $t$ such that $\mu (a) = 0$. We define the {\tmem{characteristic
polynomial}} of $a$ to be the characteristic polynomial $\chi \in \mathbb{K}
((x)) [t]$ of the multiplication endomorphism by $a$ in $\mathbb{K} ((x)) [y]
/ (f)$, so we have $\chi (a) = 0$. Since $\mathbb{K} ((x)) [y] / (f)$ is a
field, $\mu$ is irreducible and we have
\[ \chi = \mu^{n / d}, \text{ where } d \assign \deg_t \mu . \]
In addition, we know from Proposition~\ref{pp:res-charpoly} that $\chi (0) =
\pm \tmop{Res}_y (f, A)$. For the next proposition the slight abuse of
notation $\tmop{Res}_y (f, a) \assign \tmop{Res}_y (f, A)$ will be useful.

\begin{proposition}
  \label{pp:ext}Let $f \in \mathbb{K} [[x]] [y]$ be a monic irreducible
  polynomial of degree $n \geqslant 1$ such that $f (0, 0) = 0$. The map
  \begin{eqnarray*}
    v : \quad \mathbb{K} ((x)) [y] / (f) & \longrightarrow & \frac{1}{n}
    \mathbb{Z} \cup \{ \infty \}\\
    a & \longmapsto & \frac{\tmop{val}_x (\tmop{Res}_y (f, a))}{n}
  \end{eqnarray*}
  is the unique valuation on $\mathbb{K} ((x)) [y] / (f)$ that extends
  $\tmop{val}_x$.
\end{proposition}

\begin{proof}
  We first prove that $v$ does define a valuation. If $v (a) = \infty$ then
  $a$ is zero because $f$ is irreducible. We also easily verify that $v (x) =
  1$. Then, Proposition~\ref{pp:res-mul} implies that $v (ab) = v (a) + v (b)$
  holds for all $a$ and $b$ in $\mathbb{K} ((x)) [y] / (f)$. Assume that $a, b
  \neq 0$ and that $v (a) \geqslant v (b)$, or equivalently that $v (a / b)
  \geqslant 0$. Let $\mu (t) = \sum_{i = 0}^d \mu_i (x) t^i$ denote the
  minimal polynomial of $c \assign a / b$, and $\chi$ its characteristic
  polynomial. By Proposition~\ref{pp:res-charpoly}, since $f$ is monic in $y$,
  we have $\chi (0) = \pm \tmop{Res}_x (f, c)$. Since $\chi = \mu^{n / d}$ we
  deduce that $\tmop{val}_x (\mu (0)) \geqslant 0$. Since $\mu$ is monic and
  irreducible in $\mathbb{K} ((x)) [t]$, Lemma~\ref{lm:irreducible} implies
  that the Newton polygon of $\mu$ has a single edge, hence that $\tmop{val}_x
  (\mu_i) \geqslant 0$ for $i = 0, \ldots, {d}$. We conclude that the minimal
  polynomial $\mu (t - 1)$ of $a / b + 1$ has all its coefficients of
  nonnegative valuation. In particular, the constant coefficient of $\chi (t -
  1)$ has nonnegative valuation, whence $v (a / b + 1) \geqslant 0$, or
  equivalently $v (a + b) \geqslant v (b)$. Finally, we have shown that $v$ is
  a valuation.
  
  For the uniqueness, now let $v$ stand for a valuation on $\mathbb{K} ((x))
  [y] / (f)$ that extends $\tmop{val}_x$. Let us consider a non-zero element
  $c$ of $\mathbb{K} ((x)) [y] / (f)$ of positive valuation, and let $\mu$ and
  $\chi$ still denote its minimal and characteristic polynomials over
  $\mathbb{K} ((x))$. Since $\mu$ is irreducible, $\mathbb{K} ((x)) [c]$ is
  isomorphic to $\mathbb{K} ((x)) [t] / (\mu (t))$ and is a subfield of
  $\mathbb{K} ((x)) [y] / (f)$, so we may endow $\mathbb{K} ((x)) [c]$ with
  $v$ and we have
  \begin{equation}
    v (c) = v (t \tmop{mod} \mu (t)) . \label{eq:ext1}
  \end{equation}
  On the other hand, we have seen that $\mu \in \mathbb{K} [[x]] [t]$ hence
  Lemma~\ref{lm:slope} (used with $t$ instead of $y$) yields
  \begin{equation}
    v (t \tmop{mod} \mu) = \frac{\tmop{val}_x (\mu (0))}{d}, \label{eq:ext2}
  \end{equation}
  whence
  \begin{equation}
    \frac{\tmop{val}_x (\mu (0))}{d} = \frac{\tmop{val}_x (\chi (0))}{n} =
    \frac{\tmop{val}_x (\tmop{Res}_y (f, c))}{n} . \label{eq:ext3}
  \end{equation}
  The combination of~\eqref{eq:ext1}, \eqref{eq:ext2}, and~\eqref{eq:ext3}
  shows that $v (c)$ must equal $\tmop{val}_x (\tmop{Res}_y (f, c)) / n$.
\end{proof}

\subsection{Integral elements}

Let $f \in \mathbb{K} [[x]] [y]$ be a monic irreducible polynomial of degree
$n$. An element $a \in \mathbb{K} ((x)) [y] / (f)$ is said to be
{\tmem{integral}} if there exists a monic polynomial $p \in \mathbb{K} [[x]]
[t]$ such that $p (a) = 0$ holds. This is equivalent to the belonging of the
minimal and characteristic polynomials of $a$ in $\mathbb{K} [[x]] [t]$,
thanks to the Gauss lemma~{\cite[Chapter~4, Theorem~2.1]{Lang2002}}. The
following proposition gives a necessary and sufficient condition for an
element to be integral, depending on its valuation. We shall use this result
in the next section.

\begin{proposition}
  \label{pp:integral}Let $f \in \mathbb{K} [[x]] [y]$ be a monic irreducible
  polynomial of degree $n$ such that $f (0, 0) = 0$, and let $v$ denote the
  valuation of $\mathbb{K} ((x)) [y] / (f)$ that extends $\tmop{val}_{x }$.
  Then, an element~$a \in \mathbb{K} ((x)) [y] / (f)$ is integral if, and only
  if, $v (a) \geqslant 0$.
\end{proposition}

\begin{proof}
  If $a$ is integral then the constant coefficient of its characteristic
  polynomial $\chi$ has nonnegative valuation. Proposition~\ref{pp:ext} yields
  $v (a) \geqslant 0$. Conversely, if $v (a) \geqslant 0$ then the constant
  coefficient $\chi (0)$ of the characteristic polynomial $\chi \in \mathbb{K}
  ((x)) [y]$ of $a$ has nonnegative valuation. The same holds for the minimal
  polynomial $\mu \in \mathbb{K} ((x)) [t]$. Since~$\mu$ is irreducible,
  Lemma~\ref{lm:irreducible} ensures that the Newton polygon of $\mu$ has a
  single edge. Consequently $\mu$ must have all its coefficients of
  nonnegative valuation, hence $\mu \in \mathbb{K} [[x]] [t]$ and by the Gauss
  lemma $a$ is integral.
\end{proof}

\subsection{Uniformizing parameters}

In the context of Proposition~\ref{pp:ext}, the valuation group
\[ \Gamma \assign v (\mathbb{K} ((x)) [y] / (f) \setminus \{ 0 \}) \]
of $v$ is included in $\frac{1}{n} \mathbb{Z}$. So $\Gamma$ has the form
$\frac{s}{r} \mathbb{Z}$ where $r, s \in \mathbb{N}$, $r$ divides $n$, and $s$
is prime to $r$. Since $v (x) = 1$ necessarily $\Gamma$ contains $\mathbb{Z}$.
This implies that $s = 1$. The following proposition ensures that $r = n$. The
integer~$n$ is called the {\tmem{ramification index}} of $v$. In the present
context, a \tmtextit{uniformizing parameter} of $v$ is defined as an element
of minimal positive valuation.

\begin{proposition}
  \label{pp:unif}Assume that $\mathbb{K}$ is algebraically closed. Let $f \in
  \mathbb{K} [[x]] [y]$ be a monic irreducible polynomial of degree $n$ such
  that $f (0, 0) = 0$, let $v$ be the valuation on $\mathbb{K} ((x)) [y] /
  (f)$ that extends $\tmop{val}_x$, and let $\tau$ be a uniformizing parameter
  of $\mathbb{K} ((x)) [y] / (f)$. Then, $\mathbb{K} ((x)) [y] / (f)$ is
  isomorphic to $\mathbb{K} ((x)) [\tau]$ and is complete for $v$. The value
  group of $v$ is $\frac{1}{n} \mathbb{Z}$. In addition, there exist power
  series $\varphi, \psi \in \mathbb{K} [[\tau]]$ such that $x = \varphi
  (\tau)$ and $y = \psi (\tau)$ hold in $\mathbb{K} ((x)) [y] / (f)$.
\end{proposition}

\begin{proof}
  Let $r$ be as above, that is $v (\tau) = 1 / r$, and let $\mu$ denote the
  minimal polynomial of~$\tau$. From Proposition~\ref{pp:integral} we know
  that $\mu \in \mathbb{K} [[x]] [t]$. Since $\mu$ is monic and irreducible,
  Lemma~\ref{lm:irreducible} implies that its Newton polygon has a single
  edge. From Proposition~\ref{pp:ext} we know that $\tmop{val}_x (\tmop{Res}_y
  (f, \tau)) > 0$, so Proposition~\ref{pp:res-charpoly} implies that $\mu (0,
  0) = 0$. The map
  \begin{eqnarray*}
    \Psi : \quad \mathbb{K} ((x)) [t] / (\mu) & \longrightarrow & \mathbb{K}
    ((x)) [y] / (f)\\
    t & \longmapsto & \tau
  \end{eqnarray*}
  is injective by construction of $\mu$, so we may endow the field $\mathbb{K}
  ((x)) [t] / (\mu)$ with $v$, that is the unique valuation on $\mathbb{K}
  ((x)) [t] / (\mu)$ that extends $\tmop{val}_x$ according to
  Proposition~\ref{pp:ext}. By Lemma~\ref{lm:slope} applied to $\mu$, the
  slope of the edge of the Newton polygon of $\mu$ is $- 1 / r$.
  
  On the one hand we endow $\mathbb{K} ((x))$ with the ultrametric absolute
  value
  \[ | a (x) | \assign \exp (- \tmop{val}_x (a (x))), \text{ for all } a (x)
     \in \mathbb{K} ((x)), \]
  and the $\mathbb{K} ((x))$-vector space $\mathbb{K} ((x)) [t]_{< \deg \mu}$
  with the ultrametric norm
  \[ \| c (t) \|_{\infty} \assign \max_{0 \leqslant i < \deg \mu}  | c_i (x)
     |, \text{ for all } c (t) \assign \sum_{i = 0}^{\deg \mu - 1} c_i (x) t^i
     . \]
  On the other hand we consider the weights $\gamma_x \assign r$ for $x$ and
  $\gamma_t \assign 1$ for $t$ on $\mathbb{K} ((x)) [t]$ and we write
  \[ t^m \tmop{rem} \mu (t) \backassign \sum_{i = 0}^{\deg \mu - 1} u_{m, i}
     (x) t^i, \text{ where } u_{m, i} (x) \in \mathbb{K} [[x]] . \]
  For all $0 \leqslant j < \deg \mu$, we verify that
  \begin{eqnarray*}
    r \tmop{val}_x (u_{m, j} (x)) + j & = & \tmop{val} (u_{m, j} (x) t^j)\\
    & \geqslant & \min_{0 \leqslant i < \deg \mu} \tmop{val} (u_{m, i} (x)
    t^i)\\
    & = & \tmop{val} (t^m \tmop{rem} \mu (t)) \qquad \text{(by
    definition of $\tmop{val}$)}\\
    & \geqslant & \tmop{val} (t^m) \qquad \text{(by
    Lemma~\ref{lm:division})}\\
    & = & m,
  \end{eqnarray*}
  in order to obtain that
  \begin{equation}
    - \log (\| t^m \tmop{rem} \mu (t) \|_{\infty}) \geqslant \frac{m - (\deg
    \mu - 1)}{r} . \label{eq:unif-iii}
  \end{equation}
  Conversely, a lower bound for $v$ in terms of $\| \cdummy \|_{\infty}$
  follows straightforwardly from the definitions: for all $c (t) \in
  \mathbb{K} ((x)) [t]_{< \deg \mu}$ we have
  \begin{eqnarray}
    v (c (\tau)) & \geqslant & \min_{0 \leqslant i < \deg \mu} v (c_i (x)
    \tau^i) \nonumber\\
    & \geqslant & \min_{0 \leqslant i < \deg \mu} \tmop{val}_x (c_i (x))
    \enspace = \enspace - \log (\| c (t) \|_{\infty}) .  \label{eq:unif-iv}
  \end{eqnarray}
  Now let $b \in \mathbb{K} ((x)) [y] / (f)$ be non-zero and set $a \assign b
  / \tau^{rv (b)}$, so we have $v (a) = 0$, that is {$\tmop{val}_x
  (\tmop{Res}_y (f, a)) = 0$}. Let $\mu_a$ denote the minimal polynomial
  of~$a$. From Proposition~\ref{pp:integral} we know that $\mu_a \in
  \mathbb{K} [[x]] [t]$. Since $\mu_a$ is monic and irreducible,
  Lemma~\ref{lm:irreducible} implies that its Newton polygon has a single edge
  whose vertices are $(0, 0)$ and $(\deg \mu_a, 0)$. By
  Lemma~\ref{lm:hensel-bis}, since $\mathbb{K}$ is algebraically closed the
  initial of~$\mu_a$ cannot have $\geqslant 2$ distinct roots, therefore we
  can write
  \[ \mu_a (t) = (t - a_0)^{\deg \mu_a} + O (x) \]
  for some $a_0 \in \mathbb{K} \setminus \{ 0 \}$. Consequently, the minimal
  polynomial of $a - a_0$ writes $t^{\deg \mu_a} + O (x)$, so it has a single
  edge with negative slope. Propositions~\ref{pp:res-charpoly}
  and~\ref{pp:ext} imply {$v (a - a_0) > 0$}. In other words $b$ rewrites in
  the form $b = a_0 \tau^{rv (b)} + c$ with
  \[ v (c) = v (\tau^{rv (b)}  (a - a_0)) = v (b) + v (a - a_0) > v (b) . \]
  By iterating this process for $c$ we can construct unique approximate
  expansions of $b$ in the form
  \[ S_l \assign \tau^{rv (b)}  \sum_{i = 0}^l a_i \tau^i, \text{ where } a_i
     \in \mathbb{K}, \text{ such that } v (b - S_l) \geqslant v (b) + \frac{l
     + 1}{r} \]
  holds for all $l \geqslant 0$. Let $\tilde{S}_l (t) \in \mathbb{K} ((x))
  [t]_{< \deg \mu}$ represent the canonical preimage of $S_l$. The series
  $(S_l)_{l \geqslant 0}$ is a Cauchy sequence in $\mathbb{K} ((x)) [\tau]$.
  Inequality~\eqref{eq:unif-iii} implies that $(\tilde{S}_l (t))_{l \geqslant
  0}$ is a Cauchy sequence for $\| \cdummy \|_{\infty}$, so it converges to an
  element $\tilde{b} (t) \in \mathbb{K} ((x)) [t]_{< \deg \mu}$ because
  $\mathbb{K} ((x)) [t]_{< \deg \mu}$ is complete for~$\| \cdummy
  \|_{\infty}$, hence $\| \tilde{S}_l (t) - \tilde{b} (t) \|_{\infty}
  \rightarrow 0$. Inequality~\eqref{eq:unif-iv} applied to $c (t) \assign
  \tilde{S}_l (t) - \tilde{b} (t)$ implies that $(S_l)_{l \geqslant 0}$
  converges to $\tilde{b} (\tau)$. Using
  \[ v (b - \tilde{b} (\tau)) \geqslant \min (v (b - S_l), v (S_l - \tilde{b}
     (\tau))) \]
  for when $l$ tends to infinity, we deduce that $v (b - \tilde{b} (\tau)) =
  \infty$ hence that
  \begin{equation}
    b = \tilde{b} (\tau) = \tau^{rv (b)}  \sum_{i \geqslant 0} a_i \tau^i
    \label{eq:unif-i} .
  \end{equation}
  In particular, this shows that $\Psi$ is surjective, whence $\deg \mu = n$.
  Inequality~\eqref{eq:unif-iii} and equation~\eqref{eq:unif-i} further lead
  to
  \[ - \log (\| \tilde{b} (t) \|_{\infty}) \geqslant \min_{i \geqslant 0} (-
     \log (\| t^{rv (b) + i} \tmop{rem} \mu (t) \|_{\infty})) \geqslant v (b)
     - \frac{\deg \mu - 1}{r} . \]
  Combined with~\eqref{eq:unif-iv} the completeness of $\mathbb{K} ((x))
  [t]_{< \deg \mu}$ for ~$\| \cdummy \|_{\infty}$ induces the completeness of
  $\mathbb{K} ((x)) [\tau] \cong \mathbb{K} ((x)) [y] / (f)$ for $v$.
  
  Expanding $x$ and $y$ in terms of power series as we have done it for $b$
  in~\eqref{eq:unif-i} yields the existence and uniqueness of the requested
  series $\varphi$ and~$\psi$ for $x$ and $y$ in terms of $\tau$. Since
  $\tmop{val} (\varphi) = r$, still for the weights $\gamma_x \assign r$ for
  $x$ and $\gamma_t \assign 1$ for $t$, Lemma~\ref{lm:weierstrass} ensures the
  existence of $u$ invertible in $\mathbb{K} [[x, t]]$ and $\tilde{\mu} (t)
  \in \mathbb{K} [[x]] [t]$ monic such that $x - \varphi (t) = u (x, t) 
  \tilde{\mu} (t)$ and $\tmop{in} (\tilde{\mu}) = t^r - cx$ where $c \in
  \mathbb{K}$. We deduce that $\mu = \tilde{\mu}$, whence $r = n$.
\end{proof}

\begin{example}
  \label{ex:rat-puiseux}If the characteristic $p$ of $\mathbb{K}$ does not
  divide $n$, then we may write
  \[ x = \tau^n  (c_1 + c_2 \tau + c_3 \tau^2 + \cdots), \]
  and verify that
  \[ \tilde{\tau} \assign \tau \left( 1 + \frac{c_2}{c_1} \tau +
     \frac{c_3}{c_1} \tau^2 + \cdots \right)^{1 / n} \]
  is a uniformizing parameter. The parametrization $x = c_1  \tilde{\tau}^n$
  and $y = \psi (\tilde{\tau})$ is called a {\tmem{rational Puiseux
  expansion}} of $f$; see details in~{\cite{Duval1989}}.
\end{example}

\begin{example}
  Let $\mathbb{K} \assign \bar{\mathbb{F}}_2$ and $f \assign x + xy + y^2$.
  Assume that there exists a uniformizing parameter $\tau = a (x) + b (x) y$,
  where $a, b \in \mathbb{K} (x)$ such that $x = \tau^2$. Then we would have
  \[ x = a (x)^2 + b (x)^2  (x + xy), \]
  whence $b (x) = 0$, so $a (x)$ cannot exist. Consequently rational Puiseux
  expansions (defined in Example~\ref{ex:rat-puiseux}) do not exist in this
  case. In fact, $y$ is a uniformizing parameter and we have
  \[ x = \frac{y^2}{1 + y} = y^2 + y^3 + y^4 + \cdots, \]
  but we cannot extract the square root of $1 + y$.
\end{example}

\section{Places and Divisors}\label{s:divisor}

The goal of this section is to introduce the notion of place and divisor of an
algebraic plane curve. We begin by revisiting the previous results on
valuations from a geometric point of view in order to achieve properties that
do not depend on the ambient coordinates. {\tmem{From now on $\mathbb{K}$ is
assumed to be algebraically closed.}}

\subsection{Integral closures and places}

Recall that a {\tmem{discrete valuation ring}} is a principal ideal domain
that admits a unique maximal ideal. We recall that $\tmop{val}_x$ is the
unique valuation on $\mathbb{K} ((x))$ such that $\tmop{val}_x (x) = 1$.

\begin{proposition}
  \label{pp:closure-val}Let $f \in \mathbb{K} [[x]] [y]$ be an irreducible
  polynomial such that $f (0, 0) = 0$, and let $\mathbb{D} \assign \mathbb{K}
  [[x, y]] / (f)$. There exists a unique valuation $v$ on $\mathbb{D}$ having
  $\mathbb{Z}$ for valuation group. In addition, the integral closure of
  $\mathbb{D}$ in its field of fractions, written $\bar{\mathbb{D}}$, is a
  discrete valuation ring whose maximal ideal is
  \[ \mathfrak{P}= \{ a \in \tmop{Frac} (\mathbb{D}) \suchthat v (a) > 0 \} .
  \]
\end{proposition}

\begin{proof}
  If $x$ divides $f$ then $f = ux$, where $u$ is invertible in $\mathbb{K}
  [[x, y]]$. In this case there clearly exists a unique valuation $v$ on
  $\mathbb{D}$ that extends $\tmop{val}_y$. In addition, we have
  $\mathbb{D}=\mathbb{K} [[y]]$, $\tmop{Frac} (\mathbb{D}) =\mathbb{K} ((y))$,
  and $\bar{\mathbb{D}} =\mathbb{D}$, so the proposition clearly holds.
  
  Now we may assume that $f (0, y) \neq 0$. By Lemma~\ref{lm:weierstrass} we
  can write $f = ug$ where $u$ is invertible in $\mathbb{K} [[x, y]]$ and $g
  \in \mathbb{K} [[x]] [y]$ is monic and irreducible. Proposition~\ref{pp:ext}
  ensures the existence and the uniqueness of a valuation $v$ on $\mathbb{D}$
  that extends $\tmop{val}_x$. In addition, we have $\tmop{Frac} (\mathbb{D})
  =\mathbb{K} ((x)) [y] / (g)$ and Proposition~\ref{pp:integral} asserts that
  an element of $\tmop{Frac} (\mathbb{D})$ has nonnegative valuation if and
  only if it is integral. In other words, the set of integral elements
  $\bar{\mathbb{D}}$ is made of the elements of $\tmop{Frac} (\mathbb{D})$
  with nonnegative valuation, so $\bar{\mathbb{D}}$ is a ring. Let~$I$ be an
  ideal of $\bar{\mathbb{D}}$ and let $\rho$ be an element of minimal positive
  valuation in $I$. If $a \in I$ then $v (a / \rho) \geqslant 0$, whence $a
  \in (\rho)$. This proves that ideals of $\bar{\mathbb{D}}$ are principal. If
  $I$ is maximal then it is necessarily generated by a uniformizing parameter
  (an element of minimal positive valuation). Thus $\bar{\mathbb{D}}$ has a
  unique maximal ideal, which concludes the proof.
\end{proof}

Recall that since $f$ is irreducible in $\mathbb{K} [[x, y]]$ the valuation on
$\mathbb{D}$ is in one-to-one correspondence with that of $\tmop{Frac}
(\mathbb{D})$. The ideal $\mathfrak{P}$ of Proposition~\ref{pp:closure-val}
does not depend on the choice of the valuation on $\mathbb{D}$: if the
coordinates are changed linearly in $f$ then the representatives of the
elements in $\mathfrak{P}$ change accordingly. The ideal $\mathfrak{P}$ is
often called the {\tmem{place}} associated to $f$. We denote the unique
valuation on $\mathbb{D}$ in the sense of Proposition~\ref{pp:closure-val} by
$w : \mathbb{D} \twoheadrightarrow \mathbb{Z} \cup \{ \infty \}$, and call it
the {\tmem{canonical valuation}}; $w^{- 1} (1)$ is the set of uniformizing
parameters.

\subsection{Divisors}

Let $F \in \mathbb{K} [x, y, z]$ be a homogeneous irreducible polynomial
defining a curve $\mathcal{C} \assign \mathcal{V}_{\mathbb{P}} (F)$. Let
$\zeta = (\zeta_x : \zeta_y : \zeta_z) \in \mathcal{C}$. Up to a change of
coordinates we may assume in the sequel that $\zeta = (0 : 0 : 1)$. An
irreducible factor $f$ of $F (x, y, 1)$ in $\mathbb{K} [[x, y]]$ defines a
canonical valuation~$w$. If $A$ and $B$ are homogeneous polynomials in
$\mathbb{K} [x, y, z]$, then we have seen that the inequality
\[ v (A (x, y, 1) \tmop{mod} f) \leqslant v (B (x, y, 1) \tmop{mod} f) \]
is independent of the choice of the coordinates by
Proposition~\ref{pp:closure-val}.

The local factorization of $F (x, y, 1)$ in $\mathbb{K} [[x, y]]$ can be
uniquely written (up to a permutation of the factors)
\[ F (x, y, 1) = uf_1 \cdots f_r, \]
where $u$ is invertible in $\mathbb{K} [[x, y]]$ and the $f_i$ are irreducible
in $\mathbb{K} [[x, y]]$. The canonical valuation on $\mathbb{D}_i \assign
\mathbb{K} [[x, y]] / (f_i)$ is written $w_i$. For $A \in \mathbb{K} [x, y,
z]$ homogeneous and prime to $F$, the {\tmem{local divisor}} associated to $A$
at $\zeta$ is the symbolic sum
\[ \tmop{Div}_{\zeta} (A) \assign w_1 (A) \mathfrak{P}_1 + \cdots + w_r (A)
   \mathfrak{P}_r, \]
where $\mathfrak{P}_i$ is the place associated to $f_i$, for $i = 1, \ldots,
r$. The point $\zeta$ is called the {\tmem{center}} of the place
$\mathfrak{P}_i$. Since $A$ is prime to $F$ we have $\tmop{Res}_y (f_i (x, y),
A (x, y, 1)) \neq 0$, so $w_i (A) \neq \infty$. The (global) {\tmem{divisor}}
associated to~$A$ is then defined by
\[ \tmop{Div} (A) \assign \sum_{\zeta \in \mathcal{V}_{\mathbb{P}} (A, F)}
   \tmop{Div}_{\zeta} (A) . \]
This sum is finite thanks to
Proposition~\ref{pp:dimension-curve-intersection}. If $B$ is another
homogeneous polynomial prime to $F$, then we further define
\[ \tmop{Div} (A / B) \assign \tmop{Div} (A) - \tmop{Div} (B) . \]
More generally, a divisor of $\mathcal{C}$ is a finite
$\mathbb{Z}$-combination of places centered at points of $\mathcal{C}$. The
set of divisors of $\mathcal{C}$ is equipped with the following partial
ordering:
\[ \sum_{\mathfrak{P}} c_{\mathfrak{P}} \mathfrak{P} \leqslant
   \sum_{\mathfrak{P}} c'_{\mathfrak{P}} \mathfrak{P} \quad
   \Longleftrightarrow \quad \forall \mathfrak{P}, c_{\mathfrak{P}} \leqslant
   c'_{\mathfrak{P}} . \]
A divisor $D$ is said to be {\tmem{positive}} (also called {\tmem{effective}})
whenever $D \geqslant 0$. The {\tmem{degree}} of a divisor is defined by:
\[ \deg \left( \sum_{\mathfrak{P}} c_{\mathfrak{P}} \mathfrak{P} \right)
   \assign \sum_{\mathfrak{P}} c_{\mathfrak{P}} . \]

\subsection{B{\'e}zout's theorem}

The degree of the divisor $\tmop{Div} (G)$ associated to a homogeneous
polynomial $G \in \mathbb{K} [x, y, z]$ prime to $F$ only depends on the
degrees of $G$ and $F$. This is stated in the following proposition, that can
be regarded as an instance of B{\'e}zout's well known theorem for counting the
number of intersection points $\mathcal{V}_{\mathbb{P}} (G, F)$ with
multiplicities.

\begin{proposition}
  \label{pp:bezout}Let $F$ be an irreducible homogeneous polynomial in
  $\mathbb{K} [x, y, z]$. For all homogeneous polynomial $G \in \mathbb{K} [x,
  y, z]$ prime to $F$ we have $\deg (\tmop{Div} (G)) = \deg G \deg F$.
\end{proposition}

\begin{proof}
  With the same notation as in Section~\ref{s:resultant}, thanks to
  Lemma~\ref{lm:generic-coords} we may assume that~A\tmrsub{1} and~A\tmrsub{2}
  hold, so we can use Lemma~\ref{lm:generic-intersection}. Let $\zeta_x$ be a
  root of $R (x, 1)$ and let $\zeta_{y, 1}, \ldots, \zeta_{y, s}$ be the roots
  of $F (\zeta_x, y)$, so that
  \[ F (\zeta_x, y) = \prod_{i = 1}^s (y - \zeta_{y, i})^{m_i} \]
  holds with $m_i \geqslant 1$ and $m_1 + \cdots + m_s = \deg F$. Thanks to
  the usual Hensel lemma (recalled at the beginning of Section~\ref{s:hensel})
  there exist unique monic polynomials $f_i \in \mathbb{K} [[x - \zeta_x]]
  [y]$ such that
  \[ F (x, y, 1) = f_1 \cdots f_s, \]
  and $f_i (\zeta_x, y) = (y - \zeta_{y, i})^{m_i}$ and $\deg f_i = m_i$ for
  $i = 1, \ldots, s$. For all $i = 1, \ldots, s$, we write $f_{i, 1}, \ldots,
  f_{i, r_i} \in \mathbb{K} [[x - \zeta_x]] [y - \zeta_{y, i}]$ the monic
  irreducible factors of $f_i$. From Proposition~\ref{pp:res-mul} we deduce
  \begin{eqnarray*}
    \tmop{val}_{x - \zeta_x} (R (x, 1)) & = & \sum_{i = 1}^s \sum_{j =
    1}^{r_i} \tmop{val}_{x - \zeta_x} (\tmop{Res}_y (G (x, y, 1), f_{i, j} (x,
    y))) .
  \end{eqnarray*}
  Let $w_{i, j}$ denote the canonical valuation of $\mathbb{K} [[x - \zeta_x]]
  [y - \zeta_{y, j}] / (f_{i, j})$. By
  Propositions~\ref{pp:res-charpoly},~\ref{pp:ext} and~\ref{pp:unif}, we have
  \[ w_{i, j} (G (x, y, 1) \tmop{mod} f_{i, j} (x, y)) = \tmop{val}_{x -
     \zeta_x} (\tmop{Res}_y (G (x, y, 1), f_{i, j} (x, y))) . \]
  It follows that
  \[ \deg \left( \sum_{i = 1}^s \tmop{Div}_{(\zeta_x : \zeta_{y, i} : 1)} (G)
     \right) = \sum_{i = 1}^s \sum_{j = 1}^{r_i} w_{i, j} (G (x, y, 1)
     \tmop{mod} f_{i, j} (x, y)) = \tmop{val}_{x - \zeta_x} (R (x, 1)) \]
  and that
  \[ \deg_x R = \sum_{R (\zeta_x, 1) = 0} \tmop{val}_{x - \zeta_x} (R (x, 1))
     = \deg (\tmop{Div} (G)) . \]
  The conclusion follows from~A\tmrsub{2} and the usual equality $\deg R =
  \deg G \deg F$, whose proof is elementary from the definition of the
  resultant; see~{\cite[Chapter~3, Proposition~3.1.5]{CasasAlvero2019}} for
  instance.
\end{proof}

\subsection{Functions on curves}

Recall that $\mathbb{K} (\mathcal{C})$ denotes the set of the rational
functions defined on $\mathcal{C}$.

\begin{proposition}
  \label{pp:zero-div}Let $A / B \in \mathbb{K} (\mathcal{C}) \setminus \{ 0
  \}$. The divisor $\tmop{Div} (A / B)$ is zero if, and only if, $A / B \in
  \mathbb{K}$.
\end{proposition}

\begin{proof}
  If $A / B \in \mathbb{K} \setminus \{ 0 \}$ then $\tmop{Div} (A / B) = 0$
  holds by definition. Conversely, assume that $\tmop{Div} (A / B) = 0$ and
  that $A \nin \mathbb{K}$. The assumption on $\tmop{Div} (A / B)$ means that
  the valuations of $A$ and~$B$ coincide for all the valuations centered at
  points of $\mathcal{V}_{\mathbb{P}} (A, F)$. By
  Proposition~\ref{pp:dimension-curve-intersection} since $A \nin \mathbb{K}$
  there exists a point $\zeta$ in $\mathcal{V}_{\mathbb{P}} (A, F)$. Let $w$
  be the canonical valuation centered at~$\zeta$ and let $\tau$ be a
  uniformizing parameter, so $A$ and $B$ locally write $A (x, y, 1) = \sum_{i
  \geqslant m} a_i \tau^i$ and $B (x, y, 1) = \sum_{i \geqslant m} b_i \tau^i$
  with $a_m \neq 0$ and $b_m \neq 0$, thanks to Proposition~\ref{pp:unif}. It
  follows that
  \[ w \left( A - \frac{a_m}{b_m} B \right) > m = w (A) . \]
  For all the other valuations $w'$, and at other points $\zeta$, we have $w'
  \left( A - \frac{a_m}{b_m} B \right) \geqslant w' (A) = w' (B)$, by
  assumption. If $A - \frac{a_m}{b_m} B$ were prime to $F$ then we would have
  \[ \deg \left( \tmop{Div} \left( A - \frac{a_m}{b_m} B \right) \right) >
     \deg (\tmop{Div} A), \]
  which it not possible by Proposition~\ref{pp:bezout}. Thus $A \in
  \mathbb{K}$ and consequently $A / B \in \mathbb{K}$.
\end{proof}

\subsection{Adjoint divisor}\label{s:adjoint-divisor}

Let $\zeta$ be a point of the curve $\mathcal{C}=\mathcal{V}_{\mathbb{P}}
(F)$. Up to a change of coordinates we may assume that $\zeta = (0 : 0 : 1)$
and write
\[ F (x, y, 1) = u (x, y) f_1 (x, y) \cdots f_r (x, y) \]
for the irreducible factorization of $F (x, y, 1)$ in $\mathbb{K} [[x, y]]$,
with $u$ invertible. Let $w_i$, $\mathfrak{P}_i$, and $\tau_i$ represent the
canonical valuation, the place, and a uniformizing parameter of $\mathbb{D}_i
\assign \mathbb{K} [[x, y]] / (f_i)$. Let $\varphi_i (\tau_i)$ and $\psi_i
(\tau_i)$ denote the expansions of the images of $x$ and $y$ in $\tmop{Frac}
(\mathbb{D}_i)$. If $F_y \assign \frac{\partial F}{\partial y} \neq 0$ then
the {\tmem{local adjoint divisor}} at $\zeta$ of the curve $\mathcal{C}$ is
defined by
\[ \mathcal{A}_{\zeta} \assign \sum_{i = 1}^r (w_i (F_y) - w_i (\varphi_i'
   (\tau_i))) \mathfrak{P}_i . \]
Otherwise, when $F_y = 0$, we set
\[ \mathcal{A}_{\zeta} \assign \sum_{i = 1}^r (w_i (F_x) - w_i (\psi_i'
   (\tau_i))) \mathfrak{P}_i, \text{ where } F_x \assign \frac{\partial
   F}{\partial x} . \]
Note that $F_x$ and $F_y$ cannot be both identically zero because $F$ is
assumed to be absolutely irreducible. In addition $\varphi_i'$ and $\psi_i'$
cannot vanish simultaneously: if it were so then the characteristic $p$ of
$\mathbb{K}$ would be positive and $\varphi_i$ and $\psi_i$ would belong to
$\mathbb{K} [[\tau_i^p]]$, which contradicts the fact that $\tau_i$ is a
uniformizing parameter. By differentiating both sides of the equality $F
(\varphi_i (\tau_i), \psi_i (\tau_i), 1) = 0$ in $\tau_i$ we obtain
\[ F_x (\varphi_i (\tau_i), \psi_i (\tau_i), 1) \varphi_i' (\tau_i) + F_y
   (\varphi_i (\tau_i), \psi_i (\tau_i), 1) \psi_i' (\tau_i) = 0 \]
whence
\[ w_i (F_y) - w_i (\varphi_i' (\tau_i)) = w_i (F_x) - w_i (\psi_i' (\tau_i))
   \text{ for } i = 1, \ldots, r. \]
Note that $\varphi_i' \neq 0$ implies $F_y (\varphi_i (\tau_i), \psi_i
(\tau_i), 1) \neq 0$, and $\psi_i' \neq 0$ implies $F_x (\varphi_i (\tau_i),
\psi_i (\tau_i), 1) \neq 0$.

On the other hand $w_i (F_y) - w_i (\varphi_i' (\tau_i))$ is left unchanged
after dilatations of the coordinates or when $y$ is replaced by $y + \alpha x$
for all $\alpha \in \mathbb{K}$. Consequently, $\mathcal{A}_{\zeta}$ is
independent of the choice of coordinates. The (global) {\tmem{adjoint
divisor}} $\mathcal{A}$ of $\mathcal{C}$ is the sum of the
$\mathcal{A}_{\zeta}$ for $\zeta$ running over the singular points
of~$\mathcal{C}$.

\begin{remark}
  \label{rk:genus}If $F_y \neq 0$ then $w_i (\varphi_i' (\tau_i)) \geqslant
  0$, so we have $\deg \mathcal{A}_{\zeta} \leqslant \deg (\tmop{Div}_{\zeta}
  (F_y))$ and therefore $\deg \mathcal{A} \leqslant \deg (\tmop{Div} (F_y))$.
  Proposition~\ref{pp:bezout} implies that $\deg (\tmop{Div} (F_y)) = \delta
  (\delta - 1)$, where $\delta \assign \deg F$. More precisely, it is known
  that $\deg \mathcal{A} \leqslant (\delta - 1)  (\delta - 2)$ and that $\deg
  \mathcal{A}$ is even, but the proofs of these facts go beyond the scope of
  the present paper. The quantity $g \assign \frac{1}{2} ((\delta - 1) 
  (\delta - 2) - \deg \mathcal{A})$ is a nonnegative integer called the
  {\tmem{genus}} of $\mathcal{C}$.
\end{remark}

\section{The Brill--Noether method}\label{s:BNmethod}

We are now ready to present the Brill--Noether method for the computation of
Riemann--Roch spaces. We still assume that $\mathbb{K}$ is algebraically
closed and $F$ still denotes an irreducible homogeneous polynomial in
$\mathbb{K} [x, y, z]$ that defines a curve written $\mathcal{C}$. We are
given a divisor $D$ of $\mathcal{C}$ and we want to compute a
$\mathbb{K}$\mbox{-}basis of the Riemann--Roch space
\[ \mathcal{L} (D) \assign \left\{ \frac{A}{B} \in \mathbb{K} (\mathcal{C})
   \setminus \{ 0 \} \suchthat \tmop{Div} (A / B) \geqslant - D \right\} \cup
   \{ 0 \} . \]
As sketched in the introduction, the Brill--Noether method divides into two
steps. First, we look for a homogeneous polynomial $H$ prime to $F$ that
satisfies $\tmop{Div} (H) \geqslant D +\mathcal{A}$, where~$\mathcal{A}$ is
the adjoint divisor of $\mathcal{C}$ (defined in the previous section).
Second, we compute a basis $G_1, \ldots, G_{\ell}$ modulo $F$ of homogeneous
polynomials $G $ of degree $\deg H$ such that $\tmop{Div} (G ) \geqslant
\tmop{Div} (H) - D$. Essentially, in order to prove the correctness of this
method, we shall show that such a polynomial $H$ does exist and then that any
$A / B \in \mathcal{L} (D)$ can be written in the form $G  / H$ in $\mathbb{K}
(\mathcal{C})$ for some homogeneous polynomial $G$ of degree $\deg H$. The
latter problem corresponds to finding $G$ homogeneous of degree $\deg H$ such
that
\begin{equation}
  AH - BG \text{ belongs to the ideal generated by } F \text{ in }
  \mathbb{K} [x, y, z] . \label{eq:noether-rewrite}
\end{equation}

\subsection{Residue theorem}

Finding a solution to~\eqref{eq:noether-rewrite} is the purpose of the
so-called {\tmem{residue theorem}}. The proof of this theorem, given in the
following paragraphs, begins with local conditions for a polynomial $A$ to
belong to the ideal $(F, B)$. The following central lemma sheds light on the
role of the adjoint divisor in the Brill--Noether method.

\begin{lemma}
  \label{lm:noether-conditions}Let $f$ be irreducible in $\mathbb{K} [[x,
  y]]$, let $\mathbb{D} \assign \mathbb{K} [[x, y]] / (f)$, let $w$ denote the
  canonical valuation on $\mathbb{D}$, let $\tau$ be a uniformizing parameter,
  and let $\varphi$ and $\psi$ be in $\mathbb{K} [[\tau]]$ such that $x =
  \varphi (\tau)$ and $y = \psi (\tau)$ hold in $\tmop{Frac} (\mathbb{D})$.
  Then, any $c \in \tmop{Frac} (\mathbb{D})$ that satisfies
  \begin{equation}
    w (c) \geqslant w (f_y) - w (\varphi' (\tau)) \text{ if } \varphi' \neq 0
    \label{eq:noether-conditions-x}
  \end{equation}
  or
  \begin{equation}
    w (c) \geqslant w (f_x) - w (\psi' (\tau)) \text{ if } \psi' \neq 0,
    \label{eq:noether-conditions-y}
  \end{equation}
  where $f_x \assign \frac{\partial f}{\partial x}$ and $f_y \assign
  \frac{\partial f}{\partial y}$, necessarily belongs to $\mathbb{D}$.
\end{lemma}

\begin{proof}
  By the same reasoning as in Section~\ref{s:adjoint-divisor}, neither
  $\varphi' = \psi' = 0$ nor $f_x (\varphi (\tau), \psi (\tau)) = f_y (\varphi
  (\tau), \psi (\tau)) = 0$ can hold. Differentiating $f (\varphi (\tau), \psi
  (\tau)) = 0$ yields
  \begin{equation}
    f_x (\varphi (\tau), \psi (\tau)) \varphi' (\tau) + f_y (\varphi (\tau),
    \psi (\tau)) \psi' (\tau) = 0, \label{eq:noether-conditions-aux-1}
  \end{equation}
  so inequalities~\eqref{eq:noether-conditions-x}
  and~\eqref{eq:noether-conditions-y} are equivalent whenever $\varphi' \neq
  0$ and $\psi' \neq 0$. In addition, note that $\varphi' = 0$ implies $f_y
  (\varphi (\tau), \psi (\tau)) = 0$, and that $\psi' = 0$ implies $f_x
  (\varphi (\tau), \psi (\tau)) = 0$.
  
  If $x$ divides $f$ then $\mathbb{D} \assign \mathbb{K} [[y]]$, we can choose
  $\tau \assign y$, $\varphi (\tau) \assign 0$, and $\psi (\tau) \assign \tau$
  so the lemma trivially holds. The case when $y$ divides $f$ behaves
  similarly. Thanks to Lemma~\ref{lm:weierstrass}, if neither $x$ nor $y$
  divide $f$, then we may assume that $f$ belongs to $\mathbb{K} [[x]] [y]$,
  is monic in $y$ of degree $n \geqslant 1$, and satisfies $f (x, 0) \neq 0$,
  hence $\mathbb{D}=\mathbb{K} [[x]] [y] / (f)$. From
  Lemma~\ref{lm:irreducible} it follows that the Newton polygon of $f$ admits
  a single edge from $(0, m)$ to $(n, 0)$ where $m \assign \tmop{val}_x (f (x,
  0)) \geqslant 1$. Based on these preliminary remarks, the proof of the lemma
  is done by decreasing induction on
  \[ \min (w (f_x), w (f_y)) \geqslant 0. \]
  The induction ends when $\min (w (f_x), w (f_y)) = 0$. In fact if $w (f_x) =
  0$, equivalently if $m = 1$ holds, then we may take $\tau \assign y$, $\psi
  (\tau) \assign \tau$. Then, inequality~\eqref{eq:noether-conditions-y}
  simplifies to $w (c) \geqslant 0$, hence $c$ can be expressed as a power
  series in $y$, whence $c \in \mathbb{D}$. The case when $w (f_y) = 0$ is
  handled in a similar fashion.
  
  From now on we assume that $m \geqslant 2$ and $n \geqslant 2$. Up to
  swapping the variables $x$ and $y$ we may further assume that $n \leqslant
  m$. If $n < m$ then we define
  \begin{equation}
    \tilde{f} (x, z) \assign \frac{f (x, xz)}{x^n} \in \mathbb{K} [[x]] [z],
    \label{eq:blowup}
  \end{equation}
  that is irreducible and monic of degree $n$ in $z$. By
  Lemma~\ref{lm:irreducible} its Newton polygon admits a single edge from $(0,
  m - n)$ to $(n, 0)$. If $n = m$ then $\frac{f (x, xz)}{x^n}$ is also
  irreducible but the first vertex of its Newton polygon is $(0, 0)$. By the
  classical Hensel lemma (recalled at the beginning of
  Section~\ref{s:hensel}), there exists $\zeta \in \mathbb{K}$ such that
  \[ \frac{f (x, xz)}{x^n} (0, z) = (z - \zeta)^n . \]
  In this case we define
  \[ \tilde{f} (x, z) \assign \frac{f (x, x (z + \zeta))}{x^n} \in \mathbb{K}
     [[x]] [z], \]
  that is irreducible and monic of degree $n$ in $z$. Its Newton polygon has a
  single edge from $(0, \tilde{m})$ to $(n, 0)$ with $\tilde{m} > 0$. For
  convenience we set $\zeta \assign 0$ when $n < m$.
  
  The ring
  \[ \tilde{\mathbb{D}} \assign \mathbb{K} [[x, z]] / (\tilde{f}) =\mathbb{K}
     [[x]] [z] / (\tilde{f}) \]
  is a discrete valuation ring: $\tilde{w} (a (x, z)) \assign w \left( a
  \left( x, \frac{y}{x} - \zeta \right) \right)$ is the canonical valuation on
  $\tilde{\mathbb{D}}$, $\tilde{\tau} (x, z) \assign \tau \left( x,
  \frac{y}{x} - \zeta \right)$ is a uniformizing parameter, and $x = \varphi
  (\tilde{\tau})$ and
  \[ z = \frac{\psi (\tilde{\tau})}{\varphi (\tilde{\tau})} - \zeta
     \backassign \tilde{\psi} (\tilde{\tau}) \]
  hold in $\tmop{Frac} (\tilde{\mathbb{D}})$.
  
  We let $\tilde{c} (x, z) \assign c \left( x, \frac{y}{x} - \zeta \right)$.
  Let us first handle the case when $\varphi' \neq 0$. We have that
  \[ \tilde{f}_z (x, z) = \frac{f_y (x, x (z + \zeta))}{x^{n - 1}} \]
  and thus inequality~\eqref{eq:noether-conditions-x} rewrites into
  \[ \tilde{w} \left( \frac{\tilde{c} (\varphi (\tilde{\tau}), \tilde{\psi}
     (\tilde{\tau}))}{\varphi (\tilde{\tau})^{n - 1}} \right) \geqslant
     \tilde{w} (\tilde{f}_z (\varphi (\tilde{\tau}), \tilde{\psi}
     (\tilde{\tau}))) - \tilde{w} (\varphi' (\tilde{\tau})) . \]
  Let us now handle the case when $\varphi' = 0$. As previously noted, we
  necessarily have $\psi' \neq 0$, $f_y (\varphi (\tau), \psi (\tau)) = 0$,
  and $f_x (\varphi (\tau), \psi (\tau)) \neq 0$. From
  \[ \tilde{f}_x (x, z) = \frac{f_x (x, x (z + \zeta)) + (z + \zeta) f_y (x, x
     (z + \zeta))}{x^n} - n \frac{f (x, x (z + \zeta))}{x^{n + 1}}, \]
  we deduce that
  \[ \tilde{f}_x (x, z) = \frac{f_x (x, x (z + \zeta))}{x^n} \text{ holds in }
     \tmop{Frac} (\tilde{\mathbb{D}}) . \]
  Inequality~\eqref{eq:noether-conditions-y} rewrites into
  \[ \tilde{w} \left( \frac{\tilde{c} (\varphi (\tilde{\tau}), \tilde{\psi}
     (\tilde{\tau}))}{\varphi (\tilde{\tau})^n} \right) \geqslant \tilde{w}
     (\tilde{f}_x (\varphi (\tilde{\tau} ), \tilde{\psi} (\tilde{\tau}))) -
     \tilde{w} (\psi' (\tilde{\tau})), \]
  which simplifies to
  \[ \tilde{w} \left( \frac{\tilde{c} (\varphi (\tilde{\tau}), \tilde{\psi}
     (\tilde{\tau}))}{\varphi (\tilde{\tau})^{n - 1}} \right) \geqslant
     \tilde{w} (\tilde{f}_x (\varphi (\tilde{\tau} ), \tilde{\psi}
     (\tilde{\tau}))) - \tilde{w} (\tilde{\psi}' (\tilde{\tau})), \]
  by using
  \[ \tilde{w} (\tilde{\psi}' (\tilde{\tau})) = \tilde{w} \left( \frac{\psi'
     (\tilde{\tau})}{\varphi (\tilde{\tau})} - \frac{\varphi' (\tilde{\tau})
     \psi (\tilde{\tau})}{\varphi (\tilde{\tau})^2} \right) = \tilde{w} \left(
     \frac{\psi' (\tilde{\tau})}{\varphi (\tilde{\tau})} \right) = \tilde{w}
     (\psi' (\tilde{\tau})) - \tilde{w} (\varphi (\tilde{\tau})) . \]
  Since $\min (\tilde{w} (\tilde{f}_x), \tilde{w} (\tilde{f}_z)) < \min (w
  (f_x), w (f_y))$, the induction hypothesis implies that $\tilde{c} (x, z) /
  x^{n - 1}$ belongs to $\tilde{\mathbb{D}}$, so it can be written
  \[ \frac{\tilde{c} (x, z)}{x^{n - 1}} = \tilde{c}_{n - 1} (x) z^{n - 1} +
     \cdots + \tilde{c}_0 (x) + q (x, z)  \tilde{f} (x, z) \]
  with the $\tilde{c}_i (x)$ taken in $\mathbb{K} [[x]]$, and where $q (x, z)
  \in \mathbb{K} ((x)) [z]$. It follows that
  \[ \tilde{c} (x, z) = \tilde{c}_{n - 1} (x)  (xz)^{n - 1} + \cdots +
     \tilde{c}_0 (x) x^{n - 1} + x^{n - 1} q (x, z)  \tilde{f} (x, z), \]
  hence that
  \[ c (x, y) = \tilde{c}_{n - 1} (x)  (y - \zeta x)^{n - 1} + \cdots +
     \tilde{c}_0 (x) x^{n - 1} + x^{n - 1} q \left( x, \frac{y}{x} - \zeta
     \right) f (x, y) . \]
  Consequently $c$ belongs to $\mathbb{D}$.
\end{proof}

\begin{proposition}
  \label{pp:noether-conditions}Let $\zeta$ be a point of $\mathcal{C}$, and
  consider two homogeneous polynomials $A$ and $B$ that are prime to $F$. If
  $\tmop{Div}_{\zeta} (A)  \geqslant \tmop{Div}_{\zeta} (B)
  +\mathcal{A}_{\zeta}$ then $A$ belongs to the ideal generated by $F$ and $B$
  in $\mathbb{K} [x, y, z]_{\zeta}$.
\end{proposition}

\begin{proof}
  Without loss of generality, up to a change of coordinates, we may assume
  that $\zeta = (0 : 0 : 1)$ and that $F_y \neq 0$. The irreducible
  factorization of $F (x, y, 1)$ in $\mathbb{K} [[x, y]]$ can be written $uf_1
  \cdots f_r$ where $u$ is invertible and the $f_i$ are irreducible monic
  polynomials in $\mathbb{K} [[x]] [y]$, for $i = 1, \ldots, r$. Let
  $\mathbb{D}_i \assign \mathbb{K} [[x, y]] / (f_i (x, y))$, let $\tau_i$ be a
  uniformizing parameter for $\mathbb{D}_i$, let $w_i$ be its canonical
  valuation, and let $\varphi_i, \psi_i \in \mathbb{K} [[\tau_i]]$ be such
  that $x = \varphi_i (\tau_i)$ and $y = \psi_i (\tau_i)$ hold in $\tmop{Frac}
  (\mathbb{D}_i)$. We further define $\hat{f}_i \assign f_1 \cdots f_{i - 1}
  f_{i + 1} \cdots f_r$ for $i = 1, \ldots, r$.
  
  By definition, the hypothesis $\tmop{Div}_{\zeta} (A)  \geqslant
  \tmop{Div}_{\zeta} (B) +\mathcal{A}_{\zeta}$ rephrases into
  \begin{equation}
    w_i \left( \frac{A}{B} \right) \geqslant w_i (F_y) - w_i (\varphi_i'
    (\tau_i)) \text{ for } i = 1, \ldots, r. \label{eq:noether-conditions-i}
  \end{equation}
  Since $u$ is invertible in $\mathbb{D}_i$, we have $w_i (u) = 0$. And since
  $F_y (x, y, 1) \tmop{mod} f_i = u \widehat{f_i}  \frac{\partial
  f_i}{\partial y} \tmop{mod} f_i$ inequality~\eqref{eq:noether-conditions-i}
  becomes:
  \[ w_i \left( \frac{1}{\hat{f}_i}  \frac{A}{B} \right) \geqslant w_i \left(
     \frac{\partial f_i}{\partial y} \right) - w_i (\varphi_i' (\tau_i))
     \text{ for } i = 1, \ldots, r. \]
  From Lemma~\ref{lm:noether-conditions} we deduce that $\frac{1}{\hat{f}_i} 
  \frac{A (x, y, 1)}{B (x, y, 1)}$ belongs to $\mathbb{D}_i$: we write $c_i
  \in \mathbb{K} [[x]] [y]_{< \deg_y f_i}$ for its canonical representative.
  Let $c (x, y)$ stand for the canonical representative of $\frac{A (x, y,
  1)}{B (x, y, 1)}$ in $\mathbb{K} ((x)) [y] / (f_1 \cdots f_r)$, that is a
  polynomial in $\mathbb{K} ((x)) [y]_{< \deg (f_1 \cdots f_r)}$. The Chinese
  remainder formula
  \[ c (x, y) = \sum_{i = 1}^r c_i  \hat{f}_i \]
  shows that $c (x, y)$ belongs to $\mathbb{K} (x, y) \cap \mathbb{K} [[x, y]]
  =\mathbb{K} [x, y]_{(0, 0)}$ modulo $F$. In other words, $A (x, y, 1)$
  belongs to the ideal $(F (x, y, 1), B (x, y, 1))$ regarded in $\mathbb{K}
  [x, y]_{(0, 0)}$.
\end{proof}

\begin{lemma}
  \label{lm:noether}Let $F$ and $B$ be coprime homogeneous polynomials in
  $\mathbb{K} [x, y, z]$ such that $F (x, y, 0)$ and $B (x, y, 0)$ have no
  common root in $\mathbb{P}^1$. Then, $z$ is a non-zero divisor in
  $\mathbb{K} [x, y, z] / (F, B)$.
\end{lemma}

\begin{proof}
  Assume that there exists $E, U, V \in \mathbb{K} [x, y, z]$ homogeneous such
  that
  \begin{equation}
    zE (x, y, z) = U (x, y, z) F (x, y, z) + V (x, y, z) B (x, y, z) .
    \label{eq:zH}
  \end{equation}
  Substituting 0 for $z$ yields
  \[ 0 = U (x, y, 0) F (x, y, 0) + V (x, y, 0) B (x, y, 0) . \]
  Since $F (x, y, 0)$ and $B (x, y, 0)$ are homogeneous and have no common
  root in $\mathbb{P}^1$, they are coprime. Therefore there exists a
  homogeneous polynomial $W (x, y)$ such that
  \begin{eqnarray*}
    U (x, y, 0) & = & W (x, y) B (x, y, 0)\\
    V (x, y, 0) & = & - W (x, y) F (x, y, 0) .
  \end{eqnarray*}
  It follows that
  \begin{eqnarray*}
    U (x, y, z) & = & W (x, y) B (x, y, z) + O (z)\\
    V (x, y, z) & = & - W (x, y) F (x, y, z) + O (z),
  \end{eqnarray*}
  where $O (z)$ is a shorthand for homogeneous polynomials in $z\mathbb{K} [x,
  y, z]$. Plugging the latter expressions into~\eqref{eq:zH}, we deduce that
  \begin{eqnarray*}
    zE (x, y, z) & = & (W (x, y) B (x, y, z) + O (z)) F (x, y, z) + (- W (x,
    y) F (x, y, z) + O (z)) B (x, y, z)\\
    & = & O (z) F (x, y, z) + O (z) B (x, y, z) .
  \end{eqnarray*}
  This shows that $E$ belongs to $(F, B)$ and concludes the proof.
\end{proof}

\begin{proposition}
  \label{pp:269}Consider two homogeneous polynomials $A$ and $B$ prime to $F$.
  If $\tmop{Div} (A) \geqslant \tmop{Div} (B) +\mathcal{A}$ then $A$ belongs
  to the ideal $(F, B)$.
\end{proposition}

\begin{proof}
  By Proposition~\ref{pp:noether-conditions}, the polynomial $A$ belongs to
  the ideal generated by $F$ and $B$ in $\mathbb{K} [x, y, z]_{\zeta}$ for all
  $\zeta \in \mathcal{C}$. Thanks to Lemma~\ref{lm:generic-coords} the
  coordinates may be changed linearly in order to ensure that~A\tmrsub{1}
  and~A\tmrsub{2} hold for $F$ and $B$ instead of $F$ and $G$. Then
  Lemma~\ref{lm:generic-intersection} implies that $\mathcal{V}_{\mathbb{P}}
  (F, B)$ is in the affine chart $z = 1$. By
  Proposition~\ref{pp:zero-dim-ideal} the polynomial $A (x, y, 1)$ belongs to
  the ideal $(F (x, y, 1), B (x, y, 1))$. Consequently, after homogenizing, we
  obtain that the polynomial $z^m A (x, y, z)$ lies in the ideal $(F (x, y,
  z), B (x, y, z))$ for some nonnegative integer $m$. Using
  Lemma~\ref{lm:noether}, we conclude that $A (x, y, z)$ belongs to this
  ideal.
\end{proof}

Two divisors $D$ and $\tilde{D}$ of $\mathcal{C}$ are said to be
{\tmem{linearly equivalent}} if there exists a rational function $A / B \in
\mathbb{K} (\mathcal{C})$ such that $D = \tilde{D} - \tmop{Div} (A / B)$.

\begin{theorem}
  {\tmem{\label{th:residue}(Residue theorem)}} Let $D$ and $\tilde{D}$ be two
  distinct linearly equivalent divisors of~$\mathcal{C}$ and let $H$ be a
  homogeneous polynomial prime to $F$. If $\tilde{D} \geqslant 0$ and
  $\tmop{Div} (H) \geqslant D +\mathcal{A}$, then there exists a homogeneous
  polynomial G prime to $F$, of degree $\deg H$, and such that $\tmop{Div} (G
  / H) = \tilde{D} - D$.
\end{theorem}

\begin{proof}
  Let $A / B \in \mathbb{K} (\mathcal{C}) \setminus \{ 0 \}$ be such that $D =
  \tilde{D} - \tmop{Div} (A / B)$. We verify that
  \[ \tmop{Div} (AH) \geqslant \tmop{Div} (A) + D +\mathcal{A}= \tmop{Div} (B)
     + \tilde{D} +\mathcal{A}. \]
  Since $\tilde{D}$ is positive, we obtain $\tmop{Div} (AH) \geqslant
  \tmop{Div} (B) +\mathcal{A}$. By Proposition~\ref{pp:269}, the polynomial
  $AH$ belongs to the ideal $(F, B)$. Consequently, there exists a homogeneous
  polynomial~$G$ such that $BG - AH \in (F)$. It follows that $G$ is prime to
  $F$ and that $\tmop{Div} (BG) = \tmop{Div} (AH)$, whence $\tmop{Div} (G / H)
  = \tmop{Div} (A / B) = \tilde{D} - D$.
\end{proof}

\subsection{Main theorem}\label{s:main}

\begin{theorem}
  \label{tm:brillnoether}Let $\mathbb{K}$ be an algebraically closed field.
  Let $\mathcal{C}$ be an irreducible plane projective curve of equation $F =
  0$, let $\mathcal{A}$ be its adjoint divisor, and let $D$ be a divisor of
  $\mathcal{C}$. If $H$ is a non\mbox{-}zero homogeneous polynomial prime to
  $F$ that satisfies $\tmop{Div} (H) \geqslant D +\mathcal{A}$, then $H$ is a
  common denominator for all the elements of $\mathcal{L} (D)$.
\end{theorem}

\begin{proof}
  Let $A / B$ be a non-zero function in $\mathcal{L} (D)$. We will show that
  there exists a homogeneous polynomial $G$ of degree $\deg H$ such that $A /
  B$ equals $G / H$ up to a constant factor.
  
  By definition of $\mathcal{L} (D)$, the divisor $\tilde{D} \assign D +
  \tmop{Div} (A / B)$ is positive. The use of Theorem~\ref{th:residue} with
  $D$, $\tilde{D}$, and $H$ yields a homogeneous polynomial $G$ of degree
  $\deg H$ such that
  \[ \tmop{Div} (G / H) = \tilde{D} - D = \tmop{Div} (A / B) . \]
  It follows that $\tmop{Div} ((G / H) / (A / B))$ is zero, so $(G / H) / (A /
  B)$ is a constant in $\mathbb{K} (\mathcal{C})$ thanks to
  Proposition~\ref{pp:zero-div}.
\end{proof}

If a common denominator $H$ is known for $\mathcal{L} (D)$, then $G_1 / H,
\ldots, G_{\ell} / H$ is a $\mathbb{K}$-basis of $\mathcal{L} (D)$ whenever
$G_1, \ldots, G_{\ell}$ is a $\mathbb{K}$-basis of the space of homogeneous
polynomials $G \in \mathbb{K} [x, y, z] / (F)$ of degree $\deg H$ such that
$\tmop{Div} (G) \geqslant \tmop{Div} (H) - D$.

\subsection{Common denominator}

So far we have completed the presentation of the second step of the
Brill--Noether method. The first step is less intricate because it reduces to
elementary linear algebra. For this purpose it is useful to be more explicit
about the representation of divisors. For a place~$\mathfrak{P}$ of
$\mathcal{C}$ centered at $(0 : 0 : 1)$ and defined by $f (x, y) = 0$, we
assume that we are given a uniformizing parameter $\tau$ for
$\mathbb{D}=\mathbb{K} [x, y] / (f)$ and we let $w$ still denote the canonical
valuation. By Proposition~\ref{pp:unif}, there exist $\varphi, \psi \in
\mathbb{K} [[\tau]]$ such that $x = \varphi (\tau)$ and $y = \psi (\tau)$ hold
in $\tmop{Frac} (\mathbb{D})$.

Let $m$ be a nonnegative integer and let $H \in \mathbb{K} [x, y, z]$ be
homogeneous of degree $d \geqslant 0$. By definition, the condition $w (H)
\geqslant m$ rewrites into
\[ \tmop{val}_{\tau} (H (\varphi (\tau), \psi (\tau), 1)) \geqslant m, \]
that corresponds to $m$ linear equations in the coefficients of $H$. More
generally, if $D$ is a positive divisor, then the condition $\tmop{Div} (H)
\geqslant D$ corresponds to $\deg D$ linear equations in the
coefficients of~$H$.

For a divisor $D = \sum_{\mathfrak{P}} c_{\mathfrak{P}} \mathfrak{P}$ we write
$D_+ \assign \sum_{\mathfrak{P}, c_{\mathfrak{P}} > 0} c_{\mathfrak{P}}
\mathfrak{P}$ for the {\tmem{positive part}} of $D$. The condition $\tmop{Div}
(H) \geqslant D +\mathcal{A}$ is equivalent to $\tmop{Div} (H) \geqslant [D
+\mathcal{A}]_+$ that is further equivalent to a linear system with $\deg ([D
+\mathcal{A}]_+)$ equations. The number of unknowns is the number of
coefficients of $H$, that is $\binom{d + 2}{2}$. Consequently this system
admits a non-zero solution as soon as
\[ \binom{d + 2}{2} > \deg ([D +\mathcal{A}]_+) . \]
In other words, common denominators $H$ of $\mathcal{L} (D)$ do exist in
sufficiently large degrees.

\subsection{Algorithmic aspects}\label{s:algo}

The Brill--Noether method is summarized in the following algorithm for when
$\mathbb{K}$ is algebraically closed.
\medskip
\custombinding{1}{\noindent}\begin{tmparmod}{0pt}{0pt}{0em}%
  \begin{tmparsep}{0em}%
    {\tmstrong{Algorithm \tmtextup{1}}}
      \begin{descriptioncompact}
        \item[\label{algo:brillnoether}Input] An irreducible plane projective
        curve $\mathcal{C}$ defined by the equation $F = 0$, and a divisor $D$
        of $\mathcal{C}$.        
        \item[Output] A basis of $\mathcal{L} (D)$.
      \end{descriptioncompact}
      \begin{enumerate}
        \item Compute the adjoint divisor $\mathcal{A}$ of $\mathcal{C}$.
        
        \item Find a homogeneous polynomial $H$ prime to $F$ such that
        $\tmop{Div} (H) \geqslant D_+ +\mathcal{A}$.
        
        \item Compute $\tmop{Div} (H) - D$.
        
        \item Compute a basis $G_1, \ldots, G_{\ell}$ in $\mathbb{K} [x, y, z]
        / (F (x, y, z))$ of the space of the homogeneous polynomials $G$ of
        degree $\deg H$ such that $\tmop{Div} (G) \geqslant \tmop{Div} (H) -
        D$.
        
        \item Return $G_1 / H, \ldots, G_{\ell} / H$.
      \end{enumerate}
  \end{tmparsep}
\end{tmparmod}{\medskip}

The existence of $H$ in step~2 has been addressed in the preceding subsection,
so the correctness of the algorithm follows from
Theorem~\ref{tm:brillnoether}. A more precise algorithmic description of each
step goes beyond the aim of present paper. In fact a software implementation
requires suitable and efficient data structures for divisors, procedures to
operate on them, and algorithms to compute $\mathcal{A}$ and $\tmop{Div} (H)$.
References to detailed algorithms are gathered in the next section.

In practice, Riemann--Roch spaces are often needed over fields $\mathbb{K}$
that are not necessarily algebraically closed. In fact if $F \in \mathbb{K}
[x, y, z]$ is homogeneous and irreducible in $\bar{\mathbb{K}} [x, y, z]$,
then $\mathcal{A}$ is left unchanged by the Galois group of $\bar{\mathbb{K}}$
over $\mathbb{K}$ and if $D$ is also left unchanged by this Galois group, then
the linear system of step~2 can be written with coefficients in $\mathbb{K}$.
Therefore $H$ can be naturally taken in $\mathbb{K} [x, y, z]$, so $\tmop{Div}
(H)$ is also preserved by this Galois group, and the linear system in step~4
can also be written with coefficients in $\mathbb{K}$. Consequently,
Algorithm~\ref{algo:brillnoether} can compute a $\bar{\mathbb{K}}$-basis of
$\mathcal{L} (D)$ represented by elements in $\mathbb{K} (x, y, z)$.

\begin{example}
  Let us take $\mathbb{K} \assign \mathbb{F}_2$ and $F (x, y, z) \assign y^3 +
  x^3 + x^2 z$. The singular locus of $\mathcal{C} \assign
  \mathcal{V}_{\mathbb{P}} (F)$ is the singleton $\{ (0 : 0 : 1) \}$. The
  series $f (x, y) \assign F (x, y, 1) = y^3 + x^3 + x^2$ in
  $\bar{\mathbb{F}}_2 [[x, y]]$ is irreducible because, for the weight $3$ for
  $x$ and $2$ for $y$, its initial form $y^3 + x^2$ is irreducible. Let
  $\mathfrak{P}$ and $w$ denote the corresponding place and canonical
  valuation, so $w (x) = 3$ and $w (y) = 2$. The adjoint divisor of
  $\mathcal{C}$ is $\mathcal{A}= (w (F_y) - w (x) + 1) \mathfrak{P}= (4 - 3 +
  1) \mathfrak{P}= 2\mathfrak{P}$. Since $w (y^2 / x) = 1$, $\tau \assign y^2
  / x$ is a uniformizing parameter for $w$. In $\mathbb{K} (x) [y] / (f (x,
  y))$ we have $x = \tau^3 + O (\tau^4)$ and $y = \tau^2 + O (\tau^3)$.
  
  Now let $\mathfrak{P}_1$ represent the place of $\mathcal{C}$ centered at
  the smooth point $(1 : 0 : 1)$. Then $\tau_1 \assign y$ is a uniformizing
  parameter and we have $x = 1 + \tau_1^3 + O (\tau_1^4)$ locally. Let us take
  $D \assign \mathfrak{P}_1$. In order to obtain a basis of $\mathcal{L} (D)$
  we search for a homogeneous polynomial $H$ such that $\tmop{Div} (H)
  \geqslant D +\mathcal{A}$. It is easy to verify that $H = y$ is suitable and
  that $\tmop{Div} (H) = D +\mathcal{A}$. The second step of the
  Brill--Noether method consists in computing a basis of homogenous
  polynomials $G (x, y, z) = ax + by + cz$ such that $\tmop{Div} (G) \geqslant
  \tmop{Div} (H) - D =\mathcal{A}$. Since $G (x, y, 1) = c + b \tau^2 + O
  (\tau^3)$ in $\mathbb{K} (x) [y] / (f (x, y))$, we find that $c$ must be $0$
  hence that $\{ x, y \}$ is a basis of solutions for $G$. Note that $x$ and
  $y$ are linearly independent in $\mathbb{K} [x, y, z] / (F (x, y, z))$.
  Finally $\{ 1, x / y \}$ is a basis of $\mathcal{L} (D)$.
\end{example}

\section{Notes}\label{s:notes}

Good expositions of the Residue theorem stated in Theorem~\ref{th:residue} can
be found in several text books about algebraic geometry, for
instance~{\cite{FultonBook,CasasAlvero2019}}. Here the term ``residue'' comes
from the notion of {\tmem{residual sets of points}} on a curve defined
in~{\cite[p.~273]{BrillNoether1874}}, and it does not correlate with the other
classical Residue theorem for differential forms~{\cite[Chapter~4,
Corollary~4.3.3]{Stichtenoth2009}}. We conclude this paper with a list of
articles dedicated to the Brill--Noether method and its software
implementation.

\paragraph{Adjoint divisor.}The notion of adjoint divisor for plane curves was
introduced by Brill and Noether~{\cite{BrillNoether1874}} for curves with only
ordinary singularities. The extension to any curve is due to
Gorenstein~{\cite{Gorenstein1952}}. Then, other definitions have been studied
by various
authors~{\cite{Keller1974,Arbarello1985,Fulton2004,AbhyankarSathaye1974,GrecoValabrega1982}}.
Discarding a few subtleties, these definitions are essentially equivalent; see
details in~{\cite{GrecoValabrega1979,GrecoValabrega1982,Fulton2004}}. The
proof given for Lemma~\ref{lm:noether-conditions} hides the usual
desingularization method by successive blow-ups (through
equation~\ref{eq:blowup}), used in~{\cite{Fulton2004}} for instance. When the
characteristic is zero, a simpler proof based on Lagrange interpolation is
given in~{\cite[Section~3.3]{AbelardBerardiniCouvreurLecerf2022}}. It would be
interesting to adapt it to positive characteristic.

\paragraph{Algorithms.}The Brill--Noether method~{\cite{BrillNoether1874}} to
compute Riemann--Roch spaces was originally restricted to curves with only
ordinary singularities. Le Brigand and Risler extended the method to arbitrary
plane curves in~{\cite{LeBrigandRisler1988}}. Their approach led to a software
implementation by Hach{\'e}~{\cite{Hache1995,Hache1996,Hache1998}}. Other
algorithms in the vein of the Brill--Noether method have been proposed by
Huang and Ierardi~{\cite{HuangIerardi1994}} still for ordinary curves, and by
Volcheck~{\cite{Volcheck1994}} and Khuri-Makdisi~{\cite{KhuriMakdisi2007}} for
computing in the Jacobian of general curves (essentially) in characteristic
zero. Over fields of any characteristic, Campillo and Farr{\'a}n designed a
method based on Hamburger--Noether expansions~{\cite{CampilloFarran2002}}.
Implementations of Brill--Noether variants for general curves is available
within the {\tmname{Singular}} computer algebra
system~{\cite{singular-brillnoether,GreuelLossenSchulze2001}}.

More recently, faster algorithms have been designed for smooth input divisors
of nodal curves~{\cite{LeGluherSpaenlehauer2020,AbelardCouvreurLecerf2020}},
of ordinary curves~{\cite{AbelardCouvreurLecerf2022}}, and then of any curve
in characteristic zero (or positive but sufficiently
large)~{\cite{AbelardBerardiniCouvreurLecerf2022}}. Up to the present time the
design of a fast practical algorithm for any curve in any characteristic and
for any divisor is an active research topic.

The family of algorithms derived from the Brill--Noether approach is often
called ``geometric''. Finally, let us mention that another family of
algorithms for Riemann--Roch spaces, called ``arithmetic'', exists. The
current state\mbox{-}of\mbox{-}the\mbox{-}art algorithm of this family is due
to Hess~{\cite{Hess2002}}: it supports any curve and any characteristic, and
relies on integral closure computations in algebraic function fields. More
references in this direction can be found
in~{\cite{AbelardCouvreurLecerf2022}}.

\end{document}